\documentclass[11pt]{article}
\usepackage{amssymb, amsfonts, amsmath }
\usepackage{latexsym}
\providecommand{\MR}{\relax\ifhmode\unskip\space\fi MR }

\newtheorem{thm}{Theorem}[section]
\newtheorem{prop}[thm]{Proposition}
\newtheorem{lem}[thm]{Lemma}
\newtheorem{cor}[thm]{Corollary}

\newtheorem{rmk}[thm]{Remark}

\newcommand{\demo}{ \noindent {\it   Proof. }}
\newcommand{\qed}{$\Box$%\bigskip
}

\title{Orbit decidability and the conjugacy problem \\ for some extensions of groups}

\author{O.\ Bogopolski \\ \small{Institute of Mathematics of SBRAS,}\\ {\small Novosibirsk, Russia} \\ \small{ and
Universit\"{a}t Dortmund} \\ \small{Fakult\"{a}t f\"{u}r Mathematik} \\ \small{Lehrstuhl VI (Algebra)} \\
\small{Vogelpothsweg 87
D-44221 Dortmund, Germany} \\ \small{e-mail: groups@math.nsc.ru} \\ \\
\\ A.\ Martino \\ \small{Dept.\ Mat.\ Apl.\ IV, Univ.\ Pol.\ Catalunya, (Barcelona, Spain)} \\ \small{e-mail:
Armando.Martino@upc.edu}  \\
\\ \\ E.\ Ventura \\ \small{Dept.\ Mat.\ Apl.\ III, Univ.\ Pol.\ Catalunya,} \\ \small{and Centre de Recerca Matem\`{a}tica}
\\ \small{Barcelona, Catalonia} \\ \small{e-mail: enric.ventura@upc.edu} }

\begin{document}

\maketitle

\begin{abstract}
%Given a short exact sequence of groups with certain conditions, $1\to F\to G\to H\to 1$, we prove that $G$ has solvable
%conjugacy problem if and only if the corresponding action subgroup $A\leqslant Aut(F)$ is orbit decidable. Via
%Brinkmann's theorem, this provides a proof of the recent result that free-by-cyclic groups have solvable conjugacy
%problem. Several examples of free-by-free and free abelian-by-free groups with solvable conjugacy problem are also
%given. On the way, the twisted conjugacy problem is solved for (virtually) surface groups, and polycyclic groups. As an
%application, an alternative solution to the conjugacy problem for the automorphism group of the free group of rank 2 is
%given. \note{ARRANGE}
%
%On the negative side, we give an example of a group with solvable conjugacy problem but unsolvable twisted conjugacy
%problem. And we construct the first known example of a $\mathbb{Z}^4$-by-free group with unsolvable conjugacy problem.
%Also, an alternative proof for the unsolvability of the conjugacy problem in Miller's free-by-free groups is provided.

Given a short exact sequence of groups with certain conditions, $1\rightarrow F\rightarrow G\rightarrow H\rightarrow
1$, we prove that $G$ has solvable conjugacy problem if and only if the corresponding action subgroup $A\leqslant
Aut(F)$ is orbit decidable. From this, we deduce that the conjugacy problem is solvable, among others, for all groups
of the form $\mathbb{Z}^2\rtimes F_m$, $F_2\rtimes F_m$, $F_n \rtimes \mathbb{Z}$, and $\mathbb{Z}^n \rtimes_A F_m$
with virtually solvable action group $A\leqslant GL_n(\mathbb{Z})$. Also, we give an easy way of constructing groups of
the form $\mathbb{Z}^4\rtimes F_n$ and $F_3\rtimes F_n$ with unsolvable conjugacy problem. On the way, we solve the
twisted conjugacy problem for virtually surface and virtually polycyclic groups, and give an example of a group with
solvable conjugacy problem but unsolvable twisted conjugacy problem. As an application, an alternative solution to the
conjugacy problem in $Aut(F_2)$ is given.

\end{abstract}

\section{Introduction}\label{intro}

Let $G$ be a group, and $u,v\in G$. The symbol $\sim$ will be used to denote standard conjugacy in $G$ ($u\sim v$ if
there exists $x\in G$ such that $v=x^{-1}ux$). In this paper, we shall work with a twisted version, which is another
equivalence relation on $G$. Given an automorphism $\varphi \in Aut(G)$, we say that $u$ and $v$ are
\emph{$\varphi$-twisted conjugated}, denoted $u\sim_{\varphi} v$, if there exists $x\in G$ such that $v=(x\varphi
)^{-1}ux$. Of course, $\sim_{Id}$ coincides with $\sim$. Reidemeister was the first author considering the relation
$\sim_{\varphi}$ (see~\cite{Reid}), which has an important role in modern Nielsen fixed point theory. A few interesting
references can be found in~\cite{Fe}, where it is proven that the number of $\varphi$-twisted conjugacy classes in a
non-elementary word-hyperbolic group is always infinite; \cite{BMMV}, where an algorithm is given for recognizing
$\varphi$-twisted conjugacy classes in free groups; and~\cite{FT}, where the notion of twisted conjugacy separability
is analyzed.

Precisely, the recognition of twisted conjugacy classes is one of the main problems focused on in the present paper.
The twisted conjugacy problem for a group $G$ consists on finding an algorithm which, given an automorphism $\varphi
\in Aut(G)$ and two elements $u,v\in G$, decides whether $v\sim_{\varphi} u$ or not. Of course, a positive solution to
the twisted conjugacy problem automatically gives a solution to the (standard) conjugacy problem, which in turn
provides a solution to the word problem. The existence of a group $G$ with solvable word problem but unsolvable
conjugacy problem is well known (see~\cite{M1}). In this direction, one of the results here is the existence of a group
with solvable conjugacy problem, but unsolvable twisted conjugacy problem (see Corollary~\ref{c1} below).

Let us mention the real motivation and origin of the present work. Several months ago, the same three authors (together
with O. Maslakova) wrote the paper~\cite{BMMV}, where they solved the twisted conjugacy problem for finitely generated
free groups and, as a corollary, a solution to the conjugacy problem for free-by-cyclic groups was also obtained. A key
ingredient in this second result was Brinkmann's theorem, saying that there is an algorithm such that, given an
automorphism $\alpha$ of a finitely generated free group and two elements $u$ and $v$, decides whether $v$ is conjugate
to some iterated image of $u$, $v\sim u\alpha^k$ (see~\cite{Br}). At some point, Susan Hermiller brought to our
attention an article due to Miller,~\cite{M1}, where he constructed a free-by-free group with unsolvable conjugacy
problem. It turns out that the full proof in~\cite{BMMV} extends perfectly well to the bigger family of free-by-free
groups without any problem, except for a single step in the argumentation: at the point in~\cite{BMMV} where we used
Brinkmann's result, a much stronger problem arises, which we call {\em orbit decidability} (see below for definitions).
This way, orbit decidability is really the unique obstruction when extending the arguments from free-by-cyclic to
free-by-free groups.

Hence, we can deduce that a free-by-free group has solvable conjugacy problem if and only if the corresponding orbit
decidability problem is solvable. This idea is formally expressed in Theorem~\ref{main}, the main result of the present
paper. In fact, Theorem~\ref{main} works not only for free-by-free groups, but for an even bigger family of groups.
And, of course, when we restrict ourselves to the free-by-cyclic case, orbit decidability becomes exactly the
algorithmic problem already solved by Brinkmann's theorem.

More precisely, Theorem~\ref{main} talks about a short exact sequence $1\to F\to G\to H\to 1$ with some conditions on
$F$ and $H$. And states that the group $G$ has solvable conjugacy problem if and only if the action subgroup
$A_G\leqslant Aut(F)$ is \emph{orbit decidable}. Here, $A_G$ is the group of elements in $G$ acting by right
conjugation on $F$; and being orbit decidable means that, given $u,v\in F$, we can algorithmically decide whether
$v\sim u\alpha$ for some $\alpha \in A_G$. Of course, if we take $F$ to be free and $H$ to be infinite cyclic, we get
free-by-cyclic groups for $G$; and the action subgroup being orbit decidable is precisely Brinkmann's theorem (see
Subsection~\ref{free}). The conditions imposed to the groups $F$ and $H$ are the twisted conjugacy problem for $F$, and
the conjugacy problem plus a technical condition about centralizers for $H$.

In light of Theorem~\ref{main}, it becomes interesting, first, to collect groups $F$ where the twisted conjugacy
problem can be solved. And then, for every such group $F$, to study the property of orbit decidability for subgroups of
$Aut(F)$: every orbit decidable (undecidable) subgroup of $Aut(F)$ will correspond to extensions of $F$ having solvable
(unsolvable) conjugacy problem. After proving Theorem~\ref{main}, this is the main goal of the paper. We first show how
to solve the twisted conjugacy problem for surface and polycyclic groups (in fact, the solution of these problems
follows easily from the solution of the conjugacy problem for polycylic and for fundamental groups of 3-manifolds,
respectively). We also prove that the solvability of this problem goes up to finite index for these classes of groups.
Then, we concentrate on the simplest examples of groups with solvable twisted conjugacy problem, namely finitely
generated free (say $F_n$) and free abelian ($\mathbb{Z}^n$). On the positive side, we find several orbit decidable
subgroups of $Aut(F_n)$ and $Aut(\mathbb{Z}^n)=GL_n(\mathbb{Z})$, with the corresponding free-by-free and free
abelian-by-free groups with solvable conjugacy problem. Probably the most interesting result in this direction is
Proposition~\ref{ab-solv} saying that virtually solvable subgroups of $GL_n(\mathbb{Z})$ are orbit decidable; or, say
in a different way, orbit undecidable subgroups of $GL_n(\mathbb{Z})$ must contain $F_2$. On the negative side, we
establish a source of orbit undecidability, certainly related with a special role of $F_2$ inside the automorphism
group. In particular, this will give orbit undecidable subgroups in $Aut(F_n)$ for $n\geqslant 3$, which will
correspond to Miller's free-by-free groups with unsolvable conjugacy problem; and orbit undecidable subgroups in
$GL_n(\mathbb{Z})$ for $n\geqslant 4$, which will correspond to the first known examples of $\mathbb{Z}^4$-by-free
groups with unsolvable conjugacy problem.

All over the paper, $F_n$ denotes the free group of finite rank $n\geqslant 0$; and $F$, $G$ and $H$ stand for
arbitrary groups. Having Theorem~\ref{main} in mind, we will use one or another of these letters to make clear which
position in the short exact sequence the reader should think about. For example, typical groups to put in the first
position of the sequence are free groups ($F_n$), and in the third position are hyperbolic groups ($H$). We will write
morphisms as acting on the right, $x\mapsto x\varphi$. In particular, the inner automorphism given by right conjugation
by $g\in G$ is denoted $\gamma_g \colon G\to G$, $x\mapsto x\gamma_g =g^{-1}xg$. As usual, $End(G)$ denotes the monoid
of endomorphisms of $G$, $Aut(G)$ the group of automorphisms of $G$, $Inn(G)$ denotes the group of inner automorphisms,
and $Out(G)=Aut(G)/Inn(G)$. Finally, we write $[u,v]=u^{-1}v^{-1}uv$, and $C_G(u)=\{ g\in G \mid gu=ug\} \leqslant G$
for the centralizer of an element $u$.

As usual, given a group $F=\langle X\mid R\rangle$ and $m$ automorphisms $\varphi_1,\ldots ,\varphi_m \in Aut(F)$, the
free extension of $F$ by $\varphi_1,\ldots ,\varphi_m$ is the group
 $$
F\rtimes_{\varphi_1,\ldots ,\varphi_m} F_m=\langle X, t_1,\ldots ,t_m \mid R,\, t_i^{-1}xt_i=x\varphi_i \quad (x\in
X,\, i=1,\ldots ,m) \rangle.
 $$
Such a group is called \emph{$F$-by-[f.g.$\,$free]}. In particular, if $m=1$ we call it \emph{$F$-by-cyclic} and denote
by $F\rtimes_{\varphi_1} \mathbb{Z}$. It is well known that a group $G$ is $F$-by-[f.g.$\,$free] if and only if it has
a normal subgroup isomorphic to $F$, with finitely generated free quotient $G/F$ (i.e. if and only if it is the middle
term in a short exact sequence of the form $1\to F\to G\to F_n \to 1$). Above, when talking about ``free-by-free" and
``free abelian-by-free" we meant \emph{[f.g.$\,$free]-by-[f.g.$\,$free]} and
\emph{[f.g.$\,$abelian]-by-[f.g.$\,$free]}, respectively. To avoid confusions, we shall keep these last names; they
make an appropriate reference to finite generation, and they stress the fact that parenthesis are relevant (note, for
example, that surface groups are both free-by-cyclic and finitely generated, but most of them are not
[f.g.\,free]-by-cyclic).

The paper is structured as follows. In Section~\ref{preliminaries} we discuss some preliminaries concerning algorithmic
issues, setting up the notation and names used later. Section~\ref{od-cp} contains the first part of the work,
developing the relation between the conjugacy problem and the concept of orbit decidability. In
Section~\ref{more-applic} the applicability of Theorem~\ref{main} is enlarged, by finding more groups with solvable
twisted conjugacy problem (Subsection~\ref{more-tcp}), and by proving that many hyperbolic groups have small
centralizers in an algorithmic sense (Subsection~\ref{more-hyp}). Then, Section~\ref{autF2} is dedicated to solve the
conjugacy problem for $Aut(F_2)$; this problem was already known to be solvable (see~\cite{B} and~\cite{B2}), but we
present here an alternative solution to illustrate an application of the techniques developed in the paper.
Sections~\ref{+} and~\ref{-} are dedicated, respectively, to several positive and negative results, namely orbit
decidable subgroups corresponding to extensions with solvable conjugacy problem, and orbit undecidable subgroups
corresponding to extensions with unsolvable conjugacy problem. Both, in the free abelian case, and in the free case.
Finally, Section~\ref{open} is dedicated to summarize and comment on some questions and open problems.

\section{Algorithmic preliminaries}\label{preliminaries}

From the algorithmic point of view, the sentence ``let $G$ be a group" is not precise enough: the algorithmic behaviour
of $G$ may depend on how $G$ is given to us. For the purposes of this paper, we assume that a group will always be
given to us in an \emph{algorithmic} way: elements must be represented by finite objects and multiplication and
inversion must be doable in an algorithmic fashion; also, morphisms between groups, $\alpha \colon F\to G$, are to be
represented by a finite amount of information in such a way that one can algorithmically compute images of elements in
$F$.

If $G$ is finitely presented, a natural way (though not the unique one) consists in giving the group $G$ by a finite
presentation, $G=\langle x_1,\ldots ,x_n \mid r_1,\ldots ,r_m \rangle$; here, elements of $G$ are represented by words
in the $x_i$'s, multiplication and inversion are the standard ones in the free group (modulo relations), and morphisms
are given by images of generators.

\medskip

Let $G$ and $\phi \colon G\to G$ be a group and an automorphism given in an algorithmic way. Also, let $F\leqslant G$
be a subgroup. The following are interesting algorithmic problems in group theory (the first two typically known as
\emph{Dehn problems}):

\begin{itemize}
\item the \textbf{word problem} for $G$, denoted WP($G$): given an element $w\in G$, decide whether it represents the
trivial element of $G$. It is well known that there exists finitely presented groups with unsolvable word problem
(see~\cite{N} or~\cite{Boo}).
\item the \textbf{conjugacy problem} for $G$, denoted CP($G$): given two elements $u,v\in G$, decide whether they are
conjugate to each other in $G$. Clearly, a solution for the CP($G$) immediately gives a solution for the WP($G$), and
it is well known the existence of finitely presented groups with solvable word problem but unsolvable conjugacy problem
(see, for example, \cite{C} or~\cite{M1}).
\item the \textbf{$\phi$-twisted conjugacy problem} for $G$, denoted TCP$_{\phi}$($G$): given two elements
$u,v\in G$, decide whether they are $\phi$-twisted conjugate to each other in $G$ (i.e. whether $v=(g\phi)^{-1}ug$ for
some $g\in G$). Note that TCP$_{Id}$($G$) is CP($G$).
\item the \textbf{(uniform) twisted conjugacy problem} for $G$, denoted TCP($G$): given an automorphism $\phi \in Aut(G)$ and two
elements $u,v\in G$, decide whether they are $\phi$-twisted conjugate to each other in $G$. This is part of a more
general problem posted by G. Makanin in Question~10.26(a) of~\cite{Mak}. Obviously, a solution for the TCP($G$)
immediately gives a solution for the TCP$_{\phi}$($G$) for all $\phi \in Aut(G)$ (in particular, CP($G$) and WP($G$)).
In section~\ref{more-tcp}, we give an example of a finitely presented group with solvable conjugacy problem but
unsolvable twisted conjugacy problem.
\item the \textbf{membership problem} for $F$ in $G$, denoted MP($F,G$): given an element $w\in G$ decide whether it
belongs to $F$ or not. There are well-known pairs $(F,G)$ with unsolvable MP($F,G$) (see~\cite{Mih} or~\cite{M1}).
\end{itemize}

Conjugacy and twisted conjugacy problems have the ``search" variants, respectively called the \textbf{conjugacy search
problem}, CSP($G$), and the \textbf{twisted conjugacy search problem}, TCSP($G$), for $G$. They consist on additionally
finding a conjugator (or twisted-conjugator) in case it exists.

\medskip

When groups are given by finite presentations and morphisms by images of generators, the ``yes" parts of the listed
problems are always solvable. Let $G=\langle x_1,\ldots ,x_n \mid r_1,\ldots ,r_m \rangle$ and $H=\langle y_1,\ldots
,y_{n'} \mid s_1,\ldots ,s_{m'} \rangle$ be two finitely presented groups, $\phi \colon G\to H$ be a morphism (given by
$x_i \mapsto w_i(y_1,\ldots ,y_{n'})$, $i=1,\ldots ,n$), and let $F\leqslant G$ be a subgroup given by a finite set of
generators $\{ f_1 (x_1,\ldots ,x_n),\ldots ,f_t (x_1,\ldots ,x_n) \}$. We have:

\begin{itemize}
\item the \textbf{``yes" part of the word problem} for $G$, denoted WP$^+$($G$): given an element of $G$ as a (reduced)
word on the generators, $w(x_1,\ldots ,x_n)\in G$, which is known to be trivial in $G$, find an expression of $w$ as a
product of conjugates of the relations $r_i$. No matter if WP($G$) is solvable or not, WP$^+$($G$) is always solvable
by brute force: enumerate the normal closure $R=\langle\!\langle r_1,\ldots ,r_m \rangle\!\rangle$ in the free group
$\langle x_1,\ldots ,x_n \mid \,\,\rangle$, and check one by one whether its elements equal $w$ as a word after
reduction; since we know that $w=_G 1$, the process will eventually terminate. Note that, without the assumption, this
is not an algorithm because if $w\neq_G 1$ then it works forever without giving any answer.

\item the \textbf{``yes" part of the conjugacy problem} for $G$, denoted CP$^+$($G$): given two elements $u(x_1,\ldots
,x_n),\, v(x_1,\ldots ,x_n) \in G$ which are known to be conjugate to each other in $G$, find $w(x_1,\ldots ,x_n)$ such
that $w^{-1}uw=_G v$. Again, no matter if CP($G$) or even WP($G$) are solvable or not, CP$^+$($G$) is always solvable
by brute force: enumerate the elements in the free group $F=\langle x_1,\ldots ,x_n \mid \,\, \rangle$ and for each one
$w(x_1,\ldots ,x_n)$ apply WP$^+$($G$) to $v^{-1}w^{-1}uw$. We know that if $v^{-1}w^{-1}uw\neq_G 1$ this process will
never terminate. But we can start them for several $w$'s and, while running, keep opening similar processes for new
$w$'s; eventually, one of them will stop telling us which $w$ satisfies $w^{-1}uw=_G v$. As before, note that, without
the assumption, this is not an algorithm because if $u$ and $v$ are not conjugate to each other in $G$ then it works
forever without giving any answer.

\item the \textbf{``yes" part of the twisted conjugacy problem} for $G$, denoted TCP$^+$($G$): it is defined, and solved
by brute force, in the exact same way as CP$^+$($G$).

\item the \textbf{``yes" part of the membership problem} for $F$ in $G$, denoted MP$^+$($F,G$): given an element $w(x_1,
\ldots ,x_n)\in G$ known to belong to $F$, express $w$ as a word on the $f_i$'s. In a similar way as above, even if
MP($F,G$) is unsolvable, MP$^+$($F,G$) is always solvable by brute force.
\end{itemize}
Note  that, using these brute force arguments, the conjugacy problem is solvable if and only if the conjugacy search
problem is solvable. Similarly, the twisted conjugacy problem is solvable if and only if the twisted conjugacy search
problem is solvable. However, the corresponding complexities may be rather different.

\medskip

Finally, let us state few more problems, which will be considered in the present paper. We have:

\begin{itemize}
\item the \textbf{coset intersection problem} for $G$, denoted CIP($G$): given two finite sets of elements $\{ a_1,\ldots
,a_r\}$ and $\{ b_1,\ldots ,b_s\}$ in $G$, and two elements $x,y\in G$ decide whether the coset intersection $xA\cap
yB$ is empty or not, where $A=\langle a_1,\ldots ,a_r \rangle \leqslant G$ and $B=\langle b_1,\ldots ,b_s\rangle
\leqslant G$.
\end{itemize}
Let $F\unlhd G$ be a normal subgroup of $G$. By normality, two elements $u,v\in G$ conjugate to each other in $G$,
either both belong to $F$ or both outside $F$. Accordingly, CP($G$) splits into two parts, one of them relevant for our
purposes:

\begin{itemize}
\item the \textbf{conjugacy problem for $G$ restricted to $F$}, denoted CP$_F$($G$): given two elements $u,v\in F$,
decide whether they are conjugate to each other in $G$. Note that a solution to CP($G$) automatically gives a solution
to CP$_F$($G$), and that this is not, in general, the same problem as CP($F$).
\end{itemize}

\noindent Now, let $F$ be a group given in an algorithmic way, and $A$ be a subgroup of $Aut(F)$. We have:

\begin{itemize}
\item the \textbf{orbit decidability problem} for $A$, denoted OD($A$): given two elements $u,v\in F$, decide whether
there exists $\varphi \in A$ such that $u\varphi$ and $v$ are conjugate to each other in $F$. As will be seen in
Section~\ref{-}, there are finitely generated subgroups $A\leqslant Aut(F)$ with unsolvable OD($A$).
\end{itemize}

\noindent Note that orbit decidability for $A\leqslant Aut(F)$ is equivalent to the existence of an algorithm which,
for any two elements $u,v\in F$, decides whether there exists $\varphi \in A\cdot Inn(F)$ such that $u\varphi =v$. In
particular, if two subgroups $A,B\leqslant Aut(F)$ satisfy $A\cdot Inn(F)=B\cdot Inn(F)$ then OD($A$) and OD($B$) are
the same problem (in particular, $A$ is orbit decidable if and only if $B$ is orbit decidable). This means that orbit
decidability is, in fact, a property of subgroups of $Out(F)$. However, we shall keep talking about $Aut(F)$ for
notational convenience.

To finish the section, let us consider an arbitrary short exact sequence of groups,
 $$
1\longrightarrow F\stackrel{\alpha}{\longrightarrow} G\stackrel{\beta}{\longrightarrow} H\longrightarrow 1.
 $$
Such a sequence is said to be \emph{algorithmic} if: (i) the groups $F$, $G$ and $H$ and the morphisms $\alpha$ and
$\beta$ are given to us in an algorithmic way, i.e. we can effectively operate in $F$, $G$ and $H$, and compute images
under $\alpha$ and $\beta$, (ii) we can compute at least one pre-image in $G$ of any element in $H$, and (iii) we can
compute pre-images in $F$ of elements in $G$ mapping to the trivial element in $H$.

The typical example (though not the unique one) of an algorithmic short exact sequence is that given by finite
presentations and images of generators. In fact, (i) is immediate, and we can use MP$^+$($G\beta, H$) to compute
pre-images in $G$ of elements in $H$, and use MP$^+$($F\alpha, G$) to compute pre-images in $F$ of elements in $G$
mapping to $1_H$.

\section{Orbit decidability and the conjugacy problem}\label{od-cp}

Consider a short exact sequence of groups,
 $$
1\longrightarrow F\stackrel{\alpha}{\longrightarrow} G\stackrel{\beta}{\longrightarrow} H\longrightarrow 1.
 $$
Since $F\alpha$ is normal in $G$, for every $g\in G$, the conjugation $\gamma_g$ of $G$ induces an automorphism of $F$,
$x\mapsto g^{-1}xg$, which will be denoted $\varphi_g \in Aut(F)$ (note that, in general, $\varphi_g$ does not belong
to $Inn(F)$). It is clear that the set of all such automorphisms,
 $$
A_{G} =\{ \varphi_g \mid g\in G\},
 $$
is a subgroup of $Aut(F)$ containing $Inn(F)$. We shall refer to it as the \emph{action subgroup} of the given short
exact sequence.

Assuming some hypothesis on the sequence and the groups involved on it, the following theorem shows that the
solvability of the conjugacy problem for $G$ is equivalent to the orbit decidability of the action subgroup
$A_{G}\leqslant Aut(F)$.

\begin{thm}\label{main}
Let
 $$
1\longrightarrow F\stackrel{\alpha}{\longrightarrow} G\stackrel{\beta}{\longrightarrow} H\longrightarrow 1.
 $$
be an algorithmic short exact sequence of groups such that
\begin{itemize}
\item[(i)] $F$ has solvable twisted conjugacy problem,
\item[(ii)] $H$ has solvable conjugacy problem, and
\item[(iii)] for every $1\neq h\in H$, the subgroup $\langle h\rangle$ has finite index in its centralizer $C_H(h)$, and
there is an algorithm which computes a finite set of coset representatives, $z_{h,1},\ldots ,z_{h,t_h}\in H$,
 $$
C_H(h)=\langle h \rangle z_{h,1}\sqcup \cdots \sqcup \langle h \rangle z_{h,t_h}.
 $$
\end{itemize}
Then, the following are equivalent:
\begin{itemize}
\item[(a)] the conjugacy problem for $G$ is solvable,
\item[(b)] the conjugacy problem for $G$ restricted to $F$ is solvable,
\item[(c)] the action subgroup $A_G =\{ \varphi_g \mid g\in G\} \leqslant Aut(F)$ is orbit decidable.
\end{itemize}
\end{thm}

\demo As usual, we shall identify $F$ with $F\alpha \leqslant G$. By definition, $x\varphi_g =g^{-1}xg$ for every $g\in
G$ and $x\in F$. So, given two elements $x,x'\in F$, finding $g\in G$ such that $x'=g^{-1}xg$ is the same as finding
$\varphi \in A_G$ such that $x'=x\varphi$. Since $A_G=A_G\cdot Inn(F)$, this is solving OD($A_G$). Hence, $(b)$ and
$(c)$ are equivalent. It is also obvious that $(a)$ implies $(b)$. So, the relevant implication is $(b)\Rightarrow
(a)$.

Assume that $(b)$ holds, let $g,g'\in G$ be two given elements in $G$ and let us decide whether they are conjugate to
each other in $G$.

Map them to $H$. Using (ii), we can decide whether $g\beta$ and $g'\beta$ are conjugate to each other in $H$. If they
are not, then $g$ and $g'$ cannot either be conjugate to each other in $G$. Otherwise, (ii) gives us an element of $H$
conjugating $g\beta$ into $g'\beta$. Compute a pre-image $u\in G$ of this element. It satisfies $(g^u)\beta =(g\beta
)^{u\beta }=g'\beta$. Now, changing $g$ to $g^{\, u}$, we may assume that $g\beta =g'\beta$. If this is the trivial
element in $H$ (which we can decide because of (ii)) then $g$ and $g'$ lie in the image of $\alpha$, and applying (b)
we are done. Hence, we are restricted to the case $g\beta =g'\beta \neq_H 1$.

Now, compute $f\in F$ such that $g'=gf$ (this is the $\alpha$-pre-image of $g^{-1}g'$). Since $g\beta \neq_H 1$, we can
use (iii) to compute elements $z_1,\ldots ,z_t\in H$ such that $C_H(g\beta)=\langle g\beta \rangle z_1\sqcup \cdots
\sqcup \langle g\beta \rangle z_t$, and then compute a pre-image $y_i \in G$ for each $z_i$, $i=1,\ldots ,t$. Note
that, by construction, the $\beta$-images of $g$ and $y_i$ (respectively $g\beta$ and $z_i$) commute in $H$ so,
$y_i^{-1}gy_i=gp_i$ for some computable $p_i\in F$.

Since $g\beta =g'\beta$, every possible conjugator of $g$ into $g'$ must map to $C_H(g\beta)$ under $\beta$ so, it must
be of the form $g^r y_i x$ for some integer $r$, some $i\in \{ 1,\ldots ,t\}$, and some $x\in F$. Hence,
 $$
gf=g'= (x^{-1}y_i^{-1}g^{-r})g(g^r y_i x)=x^{-1}(y_i^{-1}gy_i )x=x^{-1}gp_i x.
 $$
Thus, deciding whether $g$ and $g'$ are conjugate to each other in $G$ amounts to decide whether there exists $i\in \{
1,\ldots ,t\}$ and $x\in F$ satisfying $gf=x^{-1} gp_i x$, which is equivalent to $f=(g^{-1} x^{-1} g )p_i x$ and so to
$f=(x \varphi_{g} )^{-1} p_i x$. Since $i$ takes finitely many values and the previous equation means precisely
$f\sim_{\varphi_{g}} p_i$, we can algorithmically solve this problem by hypothesis (i). This completes the proof. \qed

\begin{rmk} \label{remark}\emph{
Note that in the proof of Theorem~\ref{main} we did not use the full power of hypothesis~(i). In fact, we used a
solution to TCP$_{\phi}$($F$) only for the automorphisms in the action subgroup, $\phi \in A_G$. For specific examples,
this may be a weaker assumption than a full solution to TCP($F$). }
\end{rmk}

Theorem~\ref{main} gives us a relatively big family of groups $G$ for which the conjugacy problem reduces to its
restriction to a certain normal subgroup. Now, we point out some examples of groups $F$ and $H$ satisfying hypotheses
(i)-(iii) of the Theorem, and hence providing a family of groups $G$ for which the characterization is valid.

Suppose that $F$ is a finitely generated abelian group. For any given $\phi \in Aut(F)$ and $u,v\in F$, we have
$u\sim_{\phi} v$ if and only if $u\equiv v$ modulo $Im(\phi-Id)$. Hence, TCP($F$) reduces to solving a system of linear
equations (some over $\mathbb{Z}$ and some modulo certain integers). So it is solvable. On the other hand, Theorem~1.5
in~\cite{BMMV} states that TCP($F$) is also solvable when $F$ is finitely generated free. So, finitely generated
abelian groups, and finitely generated free groups satisfy (i).

With respect to conditions (ii) and (iii), it is well known that in a free group $H$, the centralizer of an element
$1\neq h\in H$ is cyclic and generated by the \emph{root} of $h$ (i.e., the unique non proper power $\hat{h}\in H$ such
that $h$ is a positive power of $\hat{h}$), which is computable. Clearly, then, finitely generated free groups satisfy
hypotheses (ii) and (iii).

Focusing on these families of groups, Theorem~\ref{main} is talking about solvability of the conjugacy problem for
[f.g.$\,$free]-by-[f.g.$\,$free] and [f.g.$\,$abelian]-by-[f.g.$\,$free] groups, and can be restated as follows.

\begin{thm}\label{freebyfree}
Let $F=\langle x_1, \ldots ,x_n \mid R \rangle$ be a (finitely generated) free or abelian group, and let $\varphi_1,
\ldots ,\varphi_m \in Aut(F)$. Then, the [f.g.$\,$free]-by-[f.g.$\,$free] or [f.g.$\,$abelian]-by-[f.g.$\,$free] group
 $$
G=\langle x_1, \ldots ,x_n,\, t_1, \ldots ,t_m \mid R,\, t_j^{-1}x_i t_j =x_i \varphi_j,\quad i=1,\ldots ,n ,\quad
j=1,\ldots ,m \rangle
 $$
has solvable conjugacy problem if and only if $\langle \varphi_1,\ldots ,\varphi_m \rangle \leqslant Aut(F)$ is orbit
decidable. \qed
\end{thm}

%\demo Let $H$ be the free group on $\{t_1,\ldots ,t_m \}$. The map from $G$ to $H$ sending $t_j$ to $t_j$ and killing
%the $x_i$'s is well defined, and has kernel equal to $F=\langle x_1,\ldots ,x_n \rangle \unlhd G$, which is finitely
%generated free, or finitely generated abelian. So, the short exact sequence $1\to F\to G\to H\to 1$ satisfies the
%hypothesis of Theorem~\ref{main}, and its action subgroup is $A_G=\langle \varphi_1,\ldots ,\varphi_m \rangle \cdot
%Inn(F)\leqslant Aut(F)$. \qed

In the following section, applicability of Theorem~\ref{main} will be enlarged by finding more groups with solvable
twisted conjugacy problem (see Subsection~\ref{more-tcp}), and finding more groups which satisfy condition (iii) (see
Subsection~\ref{more-hyp}).

\medskip

At this point, we want to remark that the study of the conjugacy problem in the families of
[f.g.$\,$free]-by-[f.g.$\,$free] and [f.g.$\,$free abelian]-by-[f.g.$\,$free] groups was the authors' original
motivation for developing the present research. One can interpret Theorem~\ref{freebyfree} by saying that orbit
decidability is the unique possible obstruction to the solvability of the conjugacy problem within these families of
groups. So, by finding examples of orbit decidable subgroups in $Aut(F_n)$ or $Aut(\mathbb{Z}^n)= GL_n(\mathbb{Z})$ one
is in fact solving the conjugacy problem in some [f.g.$\,$free]-by-[f.g.$\,$free] or [f.g.$\,$free
abelian]-by-[f.g.$\,$free] groups; we develop this in Section~\ref{+}.

But in the literature there are known examples of [f.g.$\,$free]-by-[f.g.$\,$free] groups with unsolvable conjugacy
problem: Miller's example of such a group (see Miller~\cite{M1}) automatically gives us an example of a finitely
generated subgroup of the automorphism group of a free group, which is orbit undecidable. We discuss this in
Section~\ref{-}, where we also find orbit undecidable subgroups in $GL_n(\mathbb{Z})$ for $n\geqslant 4$; these last
examples correspond to the first known [f.g.$\,$free abelian]-by-[f.g.$\,$free] groups with unsolvable conjugacy
problem.

\medskip

To conclude the section, we point out the following consequence of Theorem~\ref{main}.

\begin{cor}\label{same-action}
Consider two algorithmic short exact sequences
 $$
1\longrightarrow F\longrightarrow G_1\longrightarrow H_1 \longrightarrow 1
 $$
 $$
1\longrightarrow F\longrightarrow G_2\longrightarrow H_2 \longrightarrow 1
 $$
both satisfying the hypotheses of Theorem~\ref{main}, and sharing the same first term. If the two action subgroups
coincide, $A_{G_1}=A_{G_2} \leqslant Aut(F)$, then $G_1$ has solvable conjugacy problem if and only if $G_2$ has
solvable conjugacy problem. \qed
\end{cor}

Note that, in the situation of Corollary~\ref{same-action}, $A_{G_1}$ and $A_{G_2}$ can be equal, even with $G_1$ and
$G_2$ being far from isomorphic (for example, choose two very different sets of generators for $A\leqslant Aut(F)$, and
consider two extensions of $F$ by a free group with free generators acting on $F$ as the chosen automorphisms). A nice
example of this fact is the case of doubles of groups: the \emph{double of $G$ over $F\leqslant G$} is the amalgamated
product of two copies of $G$ with the corresponding $F$'s identified in the natural way.

\begin{cor}\label{doubles}
Let $1\rightarrow F\rightarrow G\rightarrow H\rightarrow 1$ be a short exact sequence satisfying the hypotheses of
Theorem~\ref{main}. Then $G$ has solvable conjugacy problem if and only if the double of $G$ over $F$ has solvable
conjugacy problem.
\end{cor}

\demo Let $1\rightarrow F'\rightarrow G'\rightarrow H'\rightarrow 1$ be another copy of the same short exact sequence,
and construct the double of $G$ over $F$, say $\underset{_{F=F'}}{G*G'}$. We have the natural short exact sequence for
it,
 $$
1\longrightarrow F=F'\longrightarrow \underset{_{F=F'}}{G*G'} \longrightarrow H*H' \longrightarrow 1,
 $$
whose action subgroup clearly coincides with $A_G \leqslant Aut(F)$. So, Corollary~\ref{same-action} gives the result.
\qed

\section{Enlarging the applicability of Theorem~\ref{main}}\label{more-applic}

In this section we shall enlarge the applicability of Theorem~\ref{main}, by solving the twisted conjugacy problem for
more groups $F$ other than free and free abelian (Subsection~\ref{more-tcp}), and by finding more groups $H$ with
conditions (ii) and (iii), other than free (Subsection~\ref{more-hyp}).

\subsection{More groups with solvable twisted conjugacy problem}\label{more-tcp}

Let us begin by proving a partial converse to Theorem~\ref{main}, in the case where $H=\mathbb{Z}$. In this case,
conditions (ii) and (iii) are obviously satisfied, and the relevant part of the statement says that, for $\varphi \in
Aut(F)$, solvability of TCP$_{\varphi^r}$($F$) for every $r\in \mathbb{Z}$ and of OD($\langle \varphi \rangle$) does
imply solvability of CP($F\rtimes_{\varphi} \mathbb{Z}$), see Remark~\ref{remark}. A weaker version of the convers is
true as well.

\begin{prop}\label{weakconvers}
Let $F$ be a finitely generated group and let $\varphi \in Aut(F)$, both given in an algorithmic way (and so
$F\rtimes_{\varphi} \mathbb{Z}$). If CP($F\rtimes_{\varphi} \mathbb{Z}$) is solvable then TCP$_{\varphi}$($F$) is
solvable, and $\langle \varphi \rangle \leqslant Aut(F)$ is orbit decidable.
\end{prop}

\demo Let us assume that $F\rtimes_{\varphi} \mathbb{Z}$ has solvable conjugacy problem. Then, $\langle \varphi \rangle
\leqslant Aut(F)$ is orbit decidable by Theorem~\ref{main}~$(a)\Rightarrow (c)$ (which uses non of the hypothesis
there) applied to the short exact sequence $1\to F\to F\rtimes_{\varphi} \mathbb{Z} \to \mathbb{Z}\to 1$.

On the other hand, let $u,v\in F$. For $x\in F$, the equality $(x\varphi)^{-1}ux=v$ holds in $F$ if and only if
$x^{-1}(tu)x=tv$ holds in $F\rtimes_{\varphi} \mathbb{Z}$. So, $u\sim_{\varphi} v$ if and only if $tu$ and $tv$ are
conjugated to each other in $F\rtimes_{\varphi} \mathbb{Z}$ by some element of $F$. By hypothesis, we know how to
decide whether $tu$ and $tv$ are conjugated to each other in $F\rtimes_{\varphi} \mathbb{Z}$ by an arbitrary element,
say $t^kx$. And,
 $$
(t^k x)^{-1}(tu)(t^k x)=x^{-1}t^{-k}tut^k x =x^{-1}t(u\varphi^k)x.
 $$
Hence, $tu$ and $tv$ are conjugated in $F\rtimes_{\varphi} \mathbb{Z}$ if and only if, for some integer $k$,
$t(u\varphi^k)$ and $tv$ are conjugated to each other in $F\rtimes_{\varphi} \mathbb{Z}$ by some element of $F$. That
is, if and only if, $u\varphi^k \sim_{\varphi} v$ for some integer $k$. But equation $u=(u\varphi)^{-1}(u\varphi )u$
tells us that $u\sim_{\varphi} u\varphi$ (see Lemma~1.7 in~\cite{BMMV}). So, solving the conjugacy problem in
$F\rtimes_{\varphi} \mathbb{Z}$ for $tu$ and $tv$, is the same as deciding whether $u\sim_{\varphi} v$ in $F$. Thus,
$F$ has solvable $\varphi$-twisted conjugacy problem. \qed

Combining Theorem~\ref{main} and this weak convers, we have the following remarkable consequence.

\begin{cor}\label{ace}
Let $F$ be a finitely generated group given in an algorithmic way. All cyclic extensions $F\rtimes_{\varphi}
\mathbb{Z}$ have solvable conjugacy problem if and only if both $F$ has solvable twisted conjugacy problem and all
cyclic subgroups of $Aut(F)$ are orbit decidable. \qed
\end{cor}

This corollary is useful in order to find more groups with solvable twisted conjugacy problem. For example, using a
recent result by J.P. Preaux, we can easily prove that surface groups have solvable twisted conjugacy problem, and that
Brinkmann's result (see~\cite{Br} or Theorem~\ref{brinkmann} in sec\-tion~\ref{free} below) is also valid for surface
groups.

\begin{thm}\label{surface}
Let $F$ be the fundamental group of a closed surface. Then, TCP($F$) is solvable and all cyclic subgroups of $Aut(F)$
are orbit decidable.
\end{thm}

\demo By Theorem~1.6 in~\cite{Prx}, $F\rtimes_{\varphi} \mathbb{Z}$ has solvable conjugacy problem for every
automorphism $\varphi \in Aut(F)$. Now, Corollary~\ref{ace} gives the result. \qed

\medskip

The same argument also provides an alternative proof of the recent result that polycyclic groups have solvable twisted
conjugacy problem (see~\cite{FT}). We remind that a group $F$ is called \emph{polycyclic} when there exists a finite
sequence of subgroups $F=K_0\vartriangleright K_1 \vartriangleright \cdots \vartriangleright K_{n-1}\vartriangleright
K_n =1$, each normal in the previous one, and with every quotient $K_i /K_{i+1}$ being cyclic (finite or infinite). If
all these quotients are infinite cyclic we say that $F$ is \emph{poly-$\mathbb{Z}$}.

\begin{thm}\label{polycylic}
Let $F$ be a polycyclic group. Then, $F$ has solvable twisted conjugacy problem, and all cyclic subgroups of $Aut(F)$
are orbit decidable.
\end{thm}

\demo Let $\varphi$ be an arbitrary automorphism of $F$. Since $F$ is normal in $F\rtimes_{\varphi} \mathbb{Z}$ with
cyclic quotient, the group $F\rtimes_{\varphi} \mathbb{Z}$ is again polycyclic. So, it has solvable conjugacy problem
(see~\cite{Rem}). Again, Corollary~\ref{ace} gives the result. \qed

\medskip

Let us study now how the twisted conjugacy problem behaves under finite index extensions. To do this, we shall need the
classical Todd-Coxeter algorithm, and a technical lemma about computing finite index characteristic subgroups.

\begin{thm}[Todd-Coxeter,~Chapter 8 in~\cite{Joh}]\label{TC}
Let $L=\langle X\mid R\,\rangle$ be a finite presentation and $K\leqslant L$ a finite index subgroup given by a finite
set of generators. There is an algorithm to compute a set of left coset representatives $\{ 1=g_0,g_1,\ldots ,g_q \}$
of $K$ in $L$ (so, $L=K\sqcup g_1 K \sqcup \cdots \sqcup g_q K$), plus the information on how to multiply them by
generators of $L$ on the left, say $x_ig_j \in g_{d(i,j)}K$. Moreover, one can algorithmically write any given $g\in L$
in the form $g=g_p k$ for some (unique) $p=0,\ldots ,q$, and some $k\in K$ expressed as a word on the generators of
$K$. In particular, MP($K,L$) is solvable.
\end{thm}

\begin{lem}\label{char-fi}
Let $F$ be a group given by a finite presentation $\langle X\mid R \,\rangle$, and suppose we are given a set of words
$\{ w_1,\ldots, w_r \}$ on $X$ such that $K=\langle w_1,\ldots ,w_r \rangle\leqslant F$ is a finite index subgroup.
Then, the finite index characteristic subgroup $K_0 =\bigcap_{[F:K']=s} K' \leqslant F$, where $s=[F:K]$, has a
computable set of generators.
\end{lem}

\demo Think $F$ as a quotient of the free group on $X$, $\pi \colon F(X)\twoheadrightarrow F$ (write $N=\ker \pi$).
Using the elementary fact that, for every finite index subgroup $L\leqslant F(X)$, $[F(X):L]\geqslant [F:L\pi ]$ with
equality if and only if $N\leqslant L$, we have
 $$
K_0 =\bigcap_{\begin{smallmatrix} K'\leqslant F \\ [F:K']= s\end{smallmatrix}} K' =\bigcap_{\begin{smallmatrix}
N\leqslant K''\leqslant F(X) \\ [F(X):K'']= s\end{smallmatrix}} (K''\pi )= \left(
\bigcap_{\begin{smallmatrix} N\leqslant K''\leqslant F(X) \\
[F(X):K'']= s\end{smallmatrix}} K'' \right) \pi.
 $$
Hence, we can compute generators for $K_0$ (as words in $X$) by first computing $s$ (use Todd-Coxeter algorithm), and
then listing and intersecting (and finally projecting into $F$) all of the finitely many subgroups $K''$ of index $s$
in $F(X)$ which contain $N$: the listing can be done by enumerating all saturated and folded Stallings graphs with $s$
vertices, and then computing a basis for the corresponding subgroup (see~\cite{St}); intersecting is easy using the
pull-back technique (again, see~\cite{St}); and deciding whether such a $K''$ contains $N$ can be done, even without
knowing explicit generators for $N$, by projecting $K''$ down to $F$, and using Todd-Coxeter again to verify if the
image has index exactly $s$, or smaller. \qed

A result by Gorjaga and Kirkinskii (see~\cite{GK}) states the existence of a group $F$ with an index two subgroup
$K\leqslant F$, such that the conjugacy problem is solvable in $K$ but unsolvable in $F$ (two years later the same was
independently proved by Collins and Miller in~\cite{CM2}, together with the opposite situation, unsolvability for $K$
and solvability for $F$). In other words, the conjugacy problem \emph{does not} go up or down through finite index
extensions. In contrast with this, the following proposition shows that the twisted conjugacy problem does (assuming
the subgroup is characteristic).

\begin{prop}\label{R-fi}
Let $F$ be a group given by a finite presentation $\langle X\mid R \rangle$, and suppose we are given a set of words
$\{ w_1,\ldots, w_r \}$ on $X$ such that $K=\langle w_1,\ldots ,w_r \rangle\leqslant F$ is a finite index subgroup.
\begin{itemize}
\item[(i)] Suppose $K$ is characteristic in $F$; if TCP($K$) is solvable then TCP($F$) is also solvable.
\item[(ii)] Suppose $K$ is normal in $F$; if TCP($K$) is solvable then CP($F$) is solvable.
\end{itemize}
\end{prop}

\demo By Todd-Coxeter algorithm, we can compute a set $\{ 1=y_1,y_2,\ldots ,y_s \}$ of left coset representatives of
$K$ in $F$, i.e. $F=K\sqcup y_2 K\sqcup \cdots \sqcup y_s K$, and use them to write any $u\in F$ (given as word in $X$)
in the form $y_i k$ for some $i=1,\ldots ,s$, and some $k\in K$ (expressed as a word in the $w_i$'s); in particular,
MP($K,F$) is solvable (see Theorem~\ref{TC}).

Suppose now that $K$ is characteristic in $F$, and TCP($K$) is solvable. Fix $\varphi \in Aut(F)$, and let $u,v\in F$
be two elements in $F$. With the previous procedure, we can write them as $u=y_i z$ and $v=y_j z'$, for some
$i,j=1,\ldots ,s$ and $z, z'\in K$. Deciding whether or not $u$ and $v$ are $\varphi$-twisted conjugated in $F$ amounts
to decide whether there exist $l=1,\ldots ,s$ and $k\in K$ such that
 $$
((ky_l^{-1})\varphi )(y_i z)(y_l k^{-1}) = y_j z'
 $$
or, equivalently,
 $$
(y_l\varphi )^{-1}y_i z y_l  = (k\varphi )^{-1} y_j z' k = y_j [(y_j^{-1}(k\varphi)y_j )^{-1} z' k].
 $$
Now observe that, if $K$ is characteristic in $F$, then the bracketed element $[(y_j^{-1}(k\varphi)y_j )^{-1} z' k]$
belongs to $K$. For every $l=1,\ldots ,s$, compute $(y_l\varphi )^{-1}y_i z y_l$ and check if it belongs to the coset
corresponding to $y_j$. If non of them do, then $y_i z$ and $y_j z'$ are not $\varphi$-twisted conjugated. Otherwise,
we have computed the non-empty list of indices $l$ and elements $k_{\, l} \in K$ such that $(y_l\varphi )^{-1}y_i z y_l
=y_j k_{\,l}$. For each one of them, it remains to decide whether there exists $k\in K$ such that $k_{\,l}=
(y_j^{-1}(k\varphi)y_j )^{-1} z' k =(k\varphi \gamma_{y_j})^{-1}z'k$. This is precisely deciding whether $k_{\,l}$ and
$z'$ are $(\varphi \gamma_{y_j})$-twisted conjugated in $K$ (which makes sense because $\varphi \gamma_{y_j}$ restricts
to an automorphism of $K$). This is doable by hypothesis, so we have proven (i).

The previous argument particularized to the special case $\varphi =Id$ works assuming only that $K$ is normal in $F$.
This shows (ii). \qed

\begin{thm}\label{groups-tcp}
Let $F$ be a group given by a finite presentation $\langle X\mid R \rangle$. If either $F$ is (finitely generated)
\begin{itemize}
\item[(i)] virtually abelian, or
\item[(ii)] virtually free, or
\item[(iii)] virtually surface group, or
\item[(iv)] virtually polycyclic,
\end{itemize}
then TCP($F$) is solvable.
\end{thm}

\demo We first need to compute generators for a finite index subgroup $K$ of $F$ being abelian, or free, or surface, or
polycyclic (and so, finitely generated, like $F$). We shall identify such a subgroup by finding a special type of
presentation for it. In the first three cases, a finite presentation will be called \emph{canonical} if the set of
relations, respectively, contains all the commutators of pairs of generators, or is empty, or consists precisely on a
single element being a surface relator. For the polycyclic case, let $L=\langle y_1,\ldots ,y_n\mid S\rangle$ be a
given finite presentation and, for every automorphism $\varphi \in Aut(L)$, consider the cyclic extension $L
\rtimes_{\varphi} \mathbb{Z}$, and its redundant presentation $\mathcal{CE}_{\varphi}(\langle y_1,\ldots ,y_n\mid
S\rangle)$ given by $L \rtimes_{\varphi} \mathbb{Z}=\langle y_1,\ldots ,y_n, t\mid S,\, t^{-1}y_i t=u_i, ty_i
t^{-1}=v_i\,\, (i=1,\ldots ,n) \rangle$, where $u_i$ and $v_i$ are words on the $y_i$'s describing, respectively, the
images and preimages of $y_i$ under $\varphi$. The point of considering such redundant presentations is that, once a
particular presentation of this type is given, one can easily verify whether $y_i\mapsto u_i$ defines an automorphism
of $G$ (with inverse given by $y_i \mapsto v_i$). Now every non-trivial, poly-$\mathbb{Z}$ group admits a presentation
of the form $\mathcal{CE}_{\varphi_n}(\cdots (\mathcal{CE}_{\varphi_1}(\langle x \mid -\rangle))\cdots )$, which we
shall call \emph{canonical} (add to this definition that the canonical presentation of the trivial group is the empty
one). Observe that, given a presentation, it can be verified whether it is canonical in any of the above four cases.

To prove the theorem, let us first enumerate the list of all subgroups of finite index in $F$, say $K_1$, $K_2,\ldots$.
This can be done by following the strategy in the proof of Lemma~\ref{char-fi}: enumerate all subgroups of the free
group $F(X)$ of a given (and increasing) index, and project them into $F$. For every $i\geqslant 1$, while computing
$K_{i+1}$, start a new parallel computation following Reidemeister-Schreier process (see Section~3.2 in \cite{MKS}) in
order to obtain a finite presentation for $K_i$, say $\langle Y_i \mid S_i \rangle$. Then, start applying to $\langle
Y_i \mid S_i \rangle$ the list of all possible sequences of Tietze transformations of any given (and increasing)
length. When one of the running processes finds a canonical presentation (of an abelian, or free, or surface, or
poly-$\mathbb{Z}$ group) then stop all of them and output this presentation. Although many of these individual
processes will never end, one of them will eventually finish because we are potentially exploring \emph{all} finite
presentations of \emph{all} finite index subgroups of $F$ and, by hypothesis, at least one of them admits a canonical
presentation (we use here the fact that every polycyclic group is virtually poly-$\mathbb{Z}$, see Proposition~2 in
Chapter~1 of~\cite{Seg}). The final output of all these parallel processes is a canonical presentation for a finite
index subgroup $K$ of $F$ being abelian, or free, or surface, or poly-$\mathbb{Z}$.

Now, apply Lemma~\ref{char-fi} to compute generators of a finite index characteristic subgroup $K_0 \leqslant F$ inside
$K$. Note that $K_0$ (for which we can obtain an explicit presentation by using Reidemeister-Schreier method) will
again be either abelian, or free, or surface, or polycyclic. So, TCP($K_0$) is solvable by results above. Hence,
TCP($F$) is also solvable by Proposition~\ref{R-fi}~(ii). \qed

\medskip

The following results are two other interesting consequences of Propositions~\ref{R-fi} and~\ref{weakconvers}.

\begin{cor}\label{c1}
There exists a finitely presented group $G$ with CP($G$) solvable, but TCP($G$) unsolvable.
\end{cor}

\demo Consider a finitely presented group $F$ with an index two subgroup $G\leqslant F$, such that CP($G$) is solvable
but CP($F$) is unsolvable (see, for example, Gorjaga-Kirkinskii~\cite{GK} or Collins-Miller~\cite{CM2}). Since index
two subgroups are normal, Proposition~\ref{R-fi}~(ii) implies that TCP($G$) must be unsolvable.\qed

\begin{cor}\label{c2}
There exists a finitely presented group $G$ and an automorphism $\varphi\in Aut(G)$ such that CP($G$) is solvable and
CP($G\rtimes_{\varphi} \mathbb{Z}$) is unsolvable. Conversely, there also exists a finitely presented group $G$ and an
automorphism $\varphi\in Aut(G)$ such that CP($G\rtimes_{\varphi} \mathbb{Z}$) is solvable and CP($G$) is unsolvable.
\end{cor}

\demo For the first assertion, consider a finitely presented group $G$ like in the previous theorem, with CP($G$)
solvable and TCP($G$) unsolvable. There must exist $\varphi \in Aut(G)$ with TCP$_{\varphi}$($G$) unsolvable. Then, by
Proposition~\ref{weakconvers}, CP($G\rtimes_{\varphi} \mathbb{Z}$) is unsolvable too.

For the second assertion, start with Collins-Miller example of an index two extension $G\leqslant F$ with unsolvable
CP($G$) but solvable CP($F$). In this construction, it is easy to see that $F=G\rtimes_{\varphi} C_2$ for an order two
automorphism $\varphi \in Aut(G)$. And solvability of CP($G\rtimes_{\varphi} C_2$) directly implies solvability of
CP($G\rtimes_{\varphi} \mathbb{Z}$) because the square of the stable letter belongs to the center of this last
group.\qed

\subsection{Centralizers in hyperbolic groups}\label{more-hyp}

Hypotheses (ii) and (iii) in Theorem~\ref{main} are also satisfied by a bigger family of groups beyond free, including
finitely generated torsion-free hyperbolic groups. All these groups $H$ provide new potential applications
of~Theorem~\ref{main}.

\begin{prop}\label{hyp}
Let $H$ be a finitely generated hyperbolic group given by a finite presentation, and let $h$ be an element of $H$.
\begin{itemize}
\item[(i)] There is an algorithm to determine whether or not the centralizer $C_H(h)$ is finite and, if it is so, to list all its
elements.
\item[(ii)] If $h$ has infinite order in $H$, then $\langle h\rangle$ has finite index in $C_H(h)$ and there is an
algorithm which computes a set of coset representatives for $\langle h\rangle$ in $C_H(h)$.
\end{itemize}
\end{prop}

\demo From the given presentation, it is possible to compute a hyperbolicity constant $\delta$ for $H$ (see~\cite{Pap}
or~\cite{HE}). Now, in~1.11 of Chapter III.$\Gamma$ of~\cite{BrH}, the authors provide an algorithm to solve the
conjugacy problem in $H$ (the reader can easily check that the q.m.c. property used there, holds with constant
$K=4\delta+1$ in our case). In fact, the same construction (with $u=v=h$) also enables us to compute a finite set of
generators for $C_H(h)$ (as labels of the closed paths at vertex $h$ in $\mathcal{G}$). So, we have computed a
generating set $\{g_1,\dots ,g_r\}$ for $C_H(h)$.

On the other hand, by~\cite{BG}, each finite subgroup of $H$ is conjugate to a subgroup which is contained in the ball
of radius $4\delta +1$ around 1, $B(4\delta +1)$. In particular, if $C_H(h)$ is finite then $C_H(h)^x \subseteq
B(4\delta +1)$ for some $x\in H$. We can apply a solution of the conjugacy search problem for $H$ (which is solvable)
to $h$ and each one of the members of $B(4\delta +1)$, ending up with a list of elements $\{ x_1,\ldots ,x_s \}$ such
that $\{ h^{x_1},\ldots ,h^{x_s} \}$ are \emph{all} the conjugates of $h$ belonging to $B(4\delta +1)$. Now, $C_H(h)$
is finite if and only if $C_H(h)^x \subseteq B(4\delta +1)$ for some $x\in H$, and this happens if and only if
$C_H(h)^{x_i} =\langle g_1^{x_i},\ldots ,g_r^{x_i} \rangle \subseteq B(4\delta +1)$ for some $i=1,\ldots ,s$ (because
if $h^y =h$ then $C_H(h)^y=C_H(h)$). Furthermore, this last condition, for a fixed $i$, is decidable in the following
way: recursively, for $l=1,2,\ldots$, check whether the set $B_l$ of all products of $l$ or less $(g_j^{x_i})^{\pm
1}$'s is contained in $B(4\delta +1)$; the process will eventually find a value of $l$ for which either $B_l \nsubseteq
B(4\delta +1)$ or $B_{l+1}=B_l\subseteq B(4\delta +1)$. In the first case $C_H(h)^{x_i} =\langle g_1^{x_i},\ldots
,g_r^{x_i} \rangle \nsubseteq B(4\delta +1)$ (and if this happens for every $i$ then $C_H(h)$ is infinite); and in the
second case $C_H(h)^{x_i} =\langle g_1^{x_i},\ldots ,g_r^{x_i} \rangle \subseteq B(4\delta +1)$, which means that
$C_H(h)$ is finite and, conjugating back by $x_i^{-1}$, we obtain the full list of its elements. This concludes the
proof of (i).

\medskip

By Corollary~3.10 in Chapter III.$\Gamma$ of ~\cite{BrH}, the group $\langle h\rangle$ has finite index in $C_H(h)$. As
explained in the proof of that corollary, different positive powers of $h$ are not conjugate to each other. Therefore,
there exists a natural number $k\leqslant |B(4\delta)|+1$, such that $h^k$ is not conjugate into the ball $B(4\delta
)$. In the proof it is claimed that each element of $C_H(h)$ lies at distance at most $2|h^k|+4\delta$ of $\langle
h^k\rangle$ and hence of $\langle h\rangle$. This means that there exists a set of coset representatives for $\langle
h\rangle$ in $C_H(h)$ inside the ball of radius $K=2(|B(4\delta)|+1)|h|+4\delta$. List the elements of such ball,
delete those not commuting with $h$ (use WP($H$) here) and, among the final list of candidates $z_1,\ldots ,z_r$, it
remains to decide which pairs $z_i,z_j$ satisfy $\langle h\rangle z_i =\langle h\rangle z_j$. This is the same as
$z=z_j z_i^{-1} \in \langle h\rangle$ which can be algorithmically checked in the following way. For every element $w$
of any group, the \textit{translation number} of $w$ is $\tau (w)=\lim_{n\to \infty} \frac{|w^n|}{n}$, where $|\cdot |$
denotes the word length with respect to a given presentation (see Definition 3.13 in Chapter III.$\Gamma$
of~\cite{BrH}). Obviously, $\tau(w)\leqslant |w|$ and, for hyperbolic groups, there exists a computable $\epsilon >0$
such that $\tau(w)>\epsilon$ for every $w$ of infinite order (see Proposition 3.1 in~\cite{Dz}, or Theorem~3.17 in
Chapter III.$\Gamma$ of~\cite{BrH}). Back to our situation, if $z=h^s$ for some integer $s$, then $\tau(z)=|s|\tau(h)$
and so, $|s|=\frac{\tau(z)}{\tau(h)}\leqslant \frac{|z|}{\epsilon}$. Thus, the exponent $s$ has finitely many
possibilities and so, we can algorithmically decide whether $z\in \langle h\rangle$. \qed

\begin{thm}\label{freebyhyp}
Consider a short exact sequence given by finite presentations,
 $$
1\longrightarrow F\longrightarrow G\longrightarrow H\longrightarrow 1,
 $$
where $F$ is (finitely generated) virtually free, or virtually abelian, or virtually surface, or virtually polycyclic
group, and $H$ is a (finitely generated) hyperbolic group where every non-trivial element of finite order has finite
centralizer. Then, $G$ has solvable conjugacy problem if and only if the corresponding action subgroup $A_G \leqslant
Aut(F)$ is orbit decidable.
\end{thm}

\demo We just need to check that hypotheses of Theorem~\ref{main} are satisfied. Certainly, (i) was already considered
above, and (ii) is well known (see~1.11 in Chapter III.$\Gamma$ of~\cite{BrH}). For (iii), apply the algorithm given in
Proposition~\ref{hyp}~(i) to $h\neq 1$. If it answers that $C_H(h)$ is finite, then the answer comes with the full list
of its elements from which, using a solution to WP($H$), we can extract the required list of coset representatives for
$C_H(h)$ modulo $\langle h\rangle$. Otherwise, $C_H(h)$ is infinite; in this case, our hypothesis ensures that $h$ has
infinite order in $H$, and hence we can apply the algorithm given in Proposition~\ref{hyp}~(ii) to find the required
list anyway. \qed

\section{The conjugacy problem for $Aut(F_2)$}\label{autF2}

The conjugacy problem for $Aut(F_n)$ is a deep open question about free groups (a possible plan to solve it has been
indicated by M. Lustig in the preprint~\cite{L}, and some partial results already published by the same author in this
direction). Among other partial results, it is known to be solvable for rank $n=2$ (see~\cite{B} or~\cite{B2}). As an
illustration of the potential applicability of Theorem~\ref{main}, we dedicate this section to deduce from it another
solution to the conjugacy problem for $Aut(F_2)$.

Consider the standard short exact sequence involving $Aut(F_2)$,
 $$
1 \to Inn(F_2) \to Aut(F_2) \to GL_2(\mathbb{Z}) \to 1.
 $$
Although $Inn(F_2)$ is isomorphic to $F_2$, which has solvable twisted conjugacy problem (see~\cite{BMMV}),
Theorem~\ref{main} cannot be directly applied to this sequence because some centralizers in $GL_2(\mathbb{Z})$ are too
large. Specifically, this group has a centre consisting of plus and minus the identity matrix. However, if we quotient
out by this order two subgroup, we obtain $PGL_2(\mathbb{Z})$ which will turn out to satisfy our required hypotheses.
Let us amend the short exact sequence above as follows.

Choose a basis $\{ a,b\}$ for $F_2$ and let $\sigma$ be the (order two) automorphism of $F_2$ sending $a$ to $a^{-1}$
and $b$ to $b^{-1}$. Note that $w\sigma =(w^{-1})^R =(w^R )^{-1}$, where $({\cdot})^R$ denotes the \textit{palindromic
reverse} of a word on $\{a,b\}$. Also, for every $\phi \in Aut(F_2)$, $\phi^R =\sigma^{-1}\phi \sigma$ acts by sending
$a$ to $(a\phi)^R$ and $b$ to $(b\phi)^R$. Consider the subgroup $\langle Inn(F_2), \sigma \rangle \leqslant Aut(F_2)$.
Since $\sigma^{-1}\phi^{-1} \sigma \phi$ abelianizes to the identity, $\phi^{-1} \sigma \phi =\sigma \gamma_x$ for some
$x\in F_2$. This (together with the elementary fact $\phi^{-1} \gamma_w \phi =\gamma_{w\phi}$) shows that $\langle
Inn(F_2), \sigma \rangle$ is normal in $Aut(F_2)$. Moreover, it is isomorphic to $F_2 \rtimes C_2$, a split extension
of $F_2$ by a cyclic group of order 2. This gives us another short exact sequence as follows
\begin{equation}\label{ses}
1 \to F_2 \rtimes C_2 \to Aut(F_2) \to PGL_2(\mathbb{Z}) \to 1.
\end{equation}
From the computational point of view we can think of~(\ref{ses}) as given by presentations of the involved groups and
the corresponding morphisms among them. But it is simpler to think of it as what is literally written: formal
expressions of the type $\sigma^r w$, where $r=0,1$ and $w\in F_2$ (in one-to-one correspondence with $\sigma^r\gamma_w
\in \langle Inn(F_2), \sigma \rangle \unlhd Aut(F_2)$); automorphisms of $F_2=\langle a,b \mid \,\,\rangle$; and
two-by-two integral matrices modulo $\pm Id$; all with the obvious morphisms among them.

Let us see that~(\ref{ses}) satisfies hypotheses (i)-(iii) of Theorem~\ref{main}. It is straightforward to show that
$F_2$ is characteristic in $F_2 \rtimes C_2$. Thus, by Proposition~\ref{R-fi}~(i) and solvability of TCP($F_2$)
(see~\cite{BMMV}), we deduce that TCP($F_2 \rtimes C_2$) is solvable. So, our short exact sequence~(\ref{ses})
satisfies~(i). On the other hand, conditions~(ii) and~(iii) follow both from standard computations with two by two
matrices, and also from considering the presentation of $PGL_2(\mathbb{Z})$ as amalgamated product
$PGL_2(\mathbb{Z})\cong D_2\underset{D_1}{\ast} D_3$ (see, for example, page 24 in~\cite{DD}).

Hence, (\ref{ses}) is a short exact sequence satisfying the hypotheses of Theorem~\ref{main}. This way, the solvability
of the conjugacy problem for $Aut(F_2)$ is equivalent to the orbit decidability of the action subgroup $A=\{
\varphi_{\phi} \mid \phi\in Aut(F_2)\} \leqslant Aut(F_2 \rtimes C_2)$. Thus, the following proposition provides a
solution to CP($Aut(F_2)$):

\begin{prop}
With the above notation, the action of $Aut(F_2)$ by right conjugation on its normal subgroup $\langle Inn(F_2), \sigma
\rangle \cong F_2 \rtimes C_2$ is orbit decidable.
\end{prop}

\demo The action of $Aut(F_2)$ on $Inn(F_2)$ is determined by the natural action on $F_2$, namely $\phi^{-1} \gamma_w
\phi =\gamma_{w\phi}$, for every $\phi\in Aut(F_2)$ and $w\in F_2$. So, by the classical Whitehead algorithm (see
Proposition~I.4.19 in~\cite{LS}), it is possible to decide whether two elements in $Inn(F_2 )$ lie in the same
$Aut(F_2)$-orbit (i.e. the action of $Aut(F_2)$ on $Inn(F_2)$ is orbit decidable). But, $Inn(F_2)$ has index 2 in
$\langle Inn(F_2), \sigma \rangle$ so, it only remains to decide whether two given elements of the form $\sigma
\gamma_u,\, \sigma \gamma_v \in \langle Inn(F_2), \sigma \rangle$, with $u,v\in F_2$, lie in the same $Aut(F_2)$-orbit
or not.

Suppose then that $\sigma\gamma _u$ and $\sigma \gamma _v$ are given to us, and let us algorithmically decide whether
there exists $\phi\in Aut(F_2)$ such that $\phi^{-1}\sigma\gamma _u\phi =\sigma \gamma _v$.

Note that, if this is the case, then $\phi^{-1}(\sigma\gamma _u )^2 \phi =(\sigma \gamma _v )^2$, while $(\sigma
\gamma_u)^2$ and $(\sigma \gamma _v)^2$ belong to $Inn(F_2)$. So, apply Whitehead algorithm to $(\sigma \gamma_u)^2$
and $(\sigma \gamma _v)^2$; if they are not in the same $Aut(F_2)$-orbit, the same is true for $\sigma \gamma_u$ and
$\sigma \gamma _v$, and we are done. Otherwise, we come up with a particular $\alpha \in Aut(F_2)$ such that
$\alpha^{-1}(\sigma\gamma _u )^2 \alpha =(\sigma \gamma _v )^2$; now, replacing $\sigma\gamma _u$ by $\alpha^{-1}
\sigma \gamma_u \alpha$, we can assume that $(\sigma \gamma_u)^2 =(\sigma \gamma _v)^2$. This means
 $$
\gamma_{(u\sigma )u} =\sigma \gamma_u \sigma \gamma_u =(\sigma \gamma_u)^2 =(\sigma \gamma_v )^2 =\sigma \gamma_v
\sigma \gamma_v =\gamma_{(v\sigma )v},
 $$
i.e. $(u\sigma )u=(v\sigma )v$. We can compute this element of $F_2$ and distinguish two cases depending whether it is
trivial or not:

\medskip

\textbf{Case 1:} $(u\sigma )u=(v\sigma )v =1$. Note that, in this case, $u$ and $v$ are palindromes (i.e. $u^R=u$ and
$v^R=v$). We shall show that $\sigma \gamma_u$ and $\sigma \gamma_v$ are always in the same $Aut(F_2)$-orbit. In fact,
notice that if $x\in \{a,b\}^{\pm 1}$ is any letter then,
 $$
\gamma_{x}^{-1} (\sigma \gamma_u )\gamma_x =\gamma_{x^{-1}} \sigma \gamma_u \gamma_{x} = \sigma \gamma_{xux}.
 $$
Hence, we can consecutively use conjugations of this type to shorten the length of $u$ (and $v$) down to 0 or 1; that
is, $\sigma \gamma_u$ (and $\sigma \gamma_v$) is in the same $Aut(F_2)$-orbit as at least one of $\sigma$, $\sigma
\gamma_a$, $\sigma \gamma_{a^{-1}}$, $\sigma \gamma_b$ or $\sigma \gamma_{b^{-1}}$. But, $\gamma_{a}^{-1}\sigma
\gamma_{a^{-1}}\gamma_{a} =\sigma \gamma_a$ and $\gamma_{b}^{-1}\sigma \gamma_{b^{-1}}\gamma_{b} =\sigma \gamma_b$.
Also, defining $\rho\in Aut(F_2)$ as $a\mapsto b$, $b\mapsto a$, we have $\rho^{-1}\sigma \gamma_a \rho =\sigma
\gamma_b$. Finally, defining $\chi\in Aut(F_2)$ as $a\mapsto ab$, $b\mapsto b$, we have $\chi^{-1} \sigma \chi =\sigma
\gamma_b$. Thus, $\sigma \gamma_u$ and $\sigma \gamma_v$ are always in the same $Aut(F_2)$-orbit in this case.

\medskip

\textbf{Case 2:} $(u\sigma )u=(v\sigma )v \neq 1$. Let $z$ be its root, i.e. $(u\sigma )u=(v\sigma )v =z^s$ for some
$s\geqslant 1$, with $z$ not being a proper power. If some $\phi \in Aut(F_2)$ satisfies $\phi^{-1} \sigma \gamma_u
\phi =\sigma \gamma_v$ then such $\phi$ must also satisfy
 $$
\gamma_{z^s \phi }=\phi^{-1} \gamma_{z^s} \phi =\phi^{-1} (\sigma \gamma_u )^2 \phi =(\sigma \gamma_v )^2 =\gamma_{z^s}
 $$
and so, $z\phi =z$. In other words, $\sigma \gamma_u$ can only be conjugated into $\sigma \gamma_v$ by automorphisms
stabilizing $z$. Using McCool's algorithm (see~ Proposition~I.5.7 in~\cite{LS}), we can compute a set of generators for
this subgroup of $Aut(F_2)$, say $Stab(z)=\langle \phi_1, \ldots,\phi_k\rangle$.

It remains to analyze how those $\phi_i$'s act on $\sigma \gamma_u$. We can compute $\phi_i^{-1}\sigma \gamma_u \phi_i$
and write this element of $F_2 \rtimes C_2$ in normal form, say $\sigma \gamma_{w_i}$. We claim that $w_i =(z\sigma
)^{n_i}u$ for some $n_i\in \mathbb{Z}$. In fact, squaring $\phi_i^{-1} \sigma \gamma_u \phi_i =\sigma \gamma_{w_i}$
(and using that $\phi_i$ stabilize $z$) we obtain
 $$
\gamma_{(u\sigma )u} =(\sigma \gamma_u )^2 =\gamma_{z^s}=\phi_i^{-1} \gamma_{z^s} \phi_i =\phi_i^{-1} (\sigma \gamma_u
)^2 \phi_i =(\sigma \gamma_{w_i} )^2 =\gamma_{(w_i \sigma )w_i}.
 $$
This implies $z^s =(u\sigma )u =(w_i \sigma )w_i$ and, applying $\sigma$ to both sides, $z^s \sigma =u(u\sigma )=w_i
(w_i \sigma )$. Now observe that
 $$
(w_i u^{-1})\gamma_{z^s \sigma} =(z^s \sigma)^{-1}(w_i u^{-1})(z^s \sigma) =(w_i^{-1}\sigma )w_i^{-1} (w_i
u^{-1})u(u\sigma )=(w_i^{-1}\sigma)(u\sigma )=w_i u^{-1},
 $$
which means that $w_i u^{-1}$ commutes with $z^s \sigma =(z\sigma )^s$. Hence, $w_i u^{-1} =(z\sigma )^{n_i}$ and $w_i
=(z\sigma )^{n_i}u$ for some computable $n_i\in \mathbb{Z}$. We have shown that, for $i=1,\ldots ,k$,
 $$
\phi_i^{-1} \sigma \gamma_u \phi_i =\sigma \gamma_{(z\sigma )^{n_i}u}.
 $$
Few more computations show that conjugating by $\phi_i^{-1}$ makes the corresponding negative effect on the exponent:
from $\phi_i^{-1} \sigma \gamma_u \phi_i =\sigma \gamma_{(z\sigma )^{n_i}u} =\gamma_{z^{n_i}}\sigma \gamma_u$ we deduce
 $$
\phi_i (\sigma \gamma_u )\phi_i^{-1} =\phi_i (\gamma_{z^{-n_i}}\phi_i^{-1} \sigma \gamma_u \phi_i
)\phi_i^{-1}=\gamma_{z^{-n_i}}\sigma \gamma_u =\sigma \gamma_{(z\sigma)^{-n_i}u}.
 $$
And conjugating by another $\phi_j$ makes an additive effect on the exponent:
 $$
\begin{array}{rcl}
\phi_j^{-1} \phi_i^{-1} \sigma \gamma_u \phi_i \phi_j & = & \phi_j^{-1} \sigma \gamma_{(z\sigma )^{n_i}u} \phi_j \\
& = & \phi_j^{-1} \gamma_{z^{n_i}} \sigma \gamma_u \phi_j \\ & = & \gamma_{z^{n_i}} \phi_j^{-1} \sigma \gamma_u \phi_j
\\ & = & \gamma_{z^{n_i}} \sigma \gamma_{(z\sigma )^{n_j}u} \\ & = & \sigma \gamma_{(z\sigma )^{n_i+n_j}u}.
\end{array}
 $$
So, conjugating by $\phi \in Stab(z)\leqslant Aut(F_2)$, we can only move from $\sigma \gamma_u$ to elements of the
form $\sigma \gamma_{(z\sigma )^{\lambda n_0}u}$, where $n_0 =\gcd (n_1,\ldots ,n_k)$, and $\lambda \in \mathbb{Z}$.
Thus, $\sigma \gamma_u$ and $\sigma \gamma_v$ belong to the same $Aut(F_2)$-orbit if and only if $v=(z\sigma )^{\lambda
n_0}u$ for some $\lambda \in \mathbb{Z}$, which happens if and only if $vu^{-1}$ is a power of $(z\sigma )^{n_0}$. This
is decidable in a free group. $\Box$

\begin{cor}
$Aut(F_2)$ has solvable conjugacy problem. $\Box$
\end{cor}

\section{Positive results}\label{+}

In this section, positive results for the free and free abelian cases are analyzed. Along the two parallel subsections,
we will give several examples of orbit decidable subgroups in $Aut(\mathbb{Z}^n)= GL_n(\mathbb{Z})$ and $Aut(F_n)$,
together with the corresponding [f.g.$\,$free abelian]-by-[f.g.$\,$free] and [f.g.$\,$free]-by-[f.g.$\,$free] groups
with solvable conjugacy problem (see Theorem~\ref{freebyfree}). The reader can easily extend all these results to
[f.g.$\,$free abelian]-by-[f.g.$\,$t.f.$\,$hyperbolic] and [f.g.$\,$free]-by-[f.g.$\,$t.f.$\,$hyperbolic] groups, by
direct application of Theorem~\ref{freebyhyp} (here, t.f. stands for torsion free).

As a technical preliminaries, we first need to solve the coset intersection problem for free and virtually free groups,
and to see that orbit decidability goes up to finite index.

\begin{prop}\label{cip-free}
Let $K$ be the free group with basis $X$. Then, the coset intersection problem CIP($K$) is solvable.
\end{prop}

\demo Without loss of generality, we can assume that $X$ is finite, since CIP($K$) only involves two finitely generated
subgroups $A,B\leqslant K$ and two words $x,y\in K$. Now, using Stallings' method (see~\cite{St}), we can construct the
corresponding core-graph $\Gamma_X(A)$ (resp. $\Gamma_X(B)$) and attach to it, and fold if necessary, a path labelled
$x$ (resp. $y$) from a new vertex $v_x$ (resp. $v_y$) to the base-point $1$ in $\Gamma_X(A)$ (resp. $1$ in
$\Gamma_X(B)$). Now, after computing the pull-back of these two finite graphs, we can easily solve the coset
intersection problem: elements from $xA\cap yB$ are precisely labels of paths from $(v_x, v_y)$ to $(1,1)$ in the
pull-back. So, $xA\cap yB \neq \emptyset$ if and only if $(v_x, v_y)$ and $(1,1)$ belong to the same connected
component of this pull-back. \qed

\begin{prop}\label{cip-fi}
Let $L$ be a group given by a finite presentation $\langle X\mid R\rangle$, and let $K$ be a finite index subgroup of
$L$ generated by a given finite set of words in $X$. If CIP($K$) is solvable then CIP($L$) is also solvable.
\end{prop}

\demo By Todd-Coxeter algorithm (see Theorem~\ref{TC}), we can compute a set of left coset representatives $\{
1=g_0,g_1,\ldots ,g_q \}$ of $K$ in $L$ (so, $L=K\sqcup g_1 K \sqcup \cdots \sqcup g_q K$), plus the information on how
to multiply them by generators of $L$ on the left, say $x_ig_j \in g_{d(i,j)}K$. This information can be used to
algorithmically write any given $g\in L$ in the form $g=g_p k$ for some (unique) $p=0,\ldots ,q$, and $k\in K$. In
particular, MP($K,L$) is solvable.

Suppose now several elements $x, y, a_1,\ldots ,a_r,$ $b_1,\ldots ,b_s\in L$ are given; we shall show how to decide
whether $xA\cap yB$ is empty or not, where $A=\langle a_1,\ldots ,a_r\rangle \leqslant L$ and $B=\langle b_1,\ldots
,b_s\rangle \leqslant L$.

Let us compute a (finite) set of left coset representatives $W$ of $A$ modulo $A\cap K$ in the following way. Enumerate
all formal reduced words in the alphabet $\{ a_1,\ldots ,a_r \}^{\pm 1}$, say $\{ w_1,w_2,\ldots \}$, starting with the
empty word, $w_1=1$, and in such a way that the length never decreases. Now, starting with the empty set $U=\emptyset$,
recursively enlarge it by adding $w_i$ whenever $w_i(A\cap K) \neq w_{i'}(A\cap K)$ for all $w_{i'} \in U$ (use here
MP($K,L$)). Since $K$ has finite index in $L$, $A\cap K$ has finite index in $A$ and the set $U$ can grow only finitely
many times. Stop the process at the moment when, for some $l$, no word of length $l$ can be added to $U$. At this
moment we have exhausted the search inside the Schreier graph of $A$ modulo $A\cap K$, and the existing list
$U=\{1=u_1,\ldots ,u_m\}$ is a set of left coset representatives of $A\cap K$ in $A$, say $A=\bigsqcup_{i=1}^{m} u_i
(A\cap K)$. Now, for every $\alpha =1,\ldots ,r$ and $i=1,\ldots ,m$, compute $d(\alpha ,i)$ such that $a_{\alpha} u_i
(A\cap K)=u_{d(\alpha ,i)}(A\cap K)$ (again using MP($K,L$)). By the Reidemeister-Schreier method (see Theorem~2.7
in~\cite{MKS}), the set $\{ u_{d(\alpha,i)}^{-1}a_{\alpha} u_i \mid \alpha =1,\ldots ,r,\,\, i=1,\ldots ,m\}$ generates
$A\cap K$. Analogously, we can compute a set $V=\{ v_1,\ldots ,v_n\}$ of left coset representatives for $B\cap K$ in
$B$, say $B=\bigsqcup_{j=1}^{n} v_j (B\cap K)$, together with a set of generators for $B\cap K$. Clearly,
 $$
xA\cap yB =\left(\bigsqcup_{i=1}^{m} xu_i(A\cap K)\right) \cap \left( \bigsqcup_{j=1}^{n} yv_j(B\cap K)\right)
=\bigsqcup_{i,j} \left( xu_i(A\cap K) \cap yv_j(B\cap K) \right).
 $$
Hence, $xA\cap yB\neq \emptyset$ is equivalent to $xu_i(A\cap K) \cap yv_j(B\cap K)\neq \emptyset$ for some $i=1,\ldots
,m$ and some $j=1,\ldots ,n$. For each $(i,j)$, consider the element $z_{i,j}=v_j^{-1}y^{-1}xu_i$ and rewrite it in the
form $g_p k$ for some (unique) $p=0,\ldots ,q$ and $k\in K$. If $p\neq 0$ (i.e. if $z_{i,j}$ does not belong to $K$)
then the intersection $xu_i(A\cap K) \cap yv_j(B\cap K)$ is empty. Otherwise, $z_{i,j}\in K$ and we are reduced to
verify whether $z_{i,j}(A\cap K) \cap (B\cap K)\neq \emptyset$. This can be done using CIP($K$). \qed

\medskip

Let us apply this result to $GL_2(\mathbb{Z})$. It is well-known that $GL_2(\mathbb{Z})$ admits the following
presentation
 \begin{equation}\label{presGL2}
GL_2(\mathbb{Z})\cong D_4\underset{D_2}{\ast} D_6 =\langle t_4,\, x_4 \mid t_4^2,\, x_4^4,\, (t_4x_4)^2 \rangle
\underset{\left\langle
\begin{array}{c} t_4=t_6 \\ x_4^2=x_6^3 \end{array}\right\rangle}{\ast} \langle t_6,\, x_6 \mid t_6^2,\, x_6^6,\,
(t_6x_6)^2 \rangle,
 \end{equation}
where $t_4=t_6=\left( \begin{smallmatrix} 0 & 1 \\ 1 & 0
\end{smallmatrix}\right)$, $x_4=\left( \begin{smallmatrix} 0 & -1 \\
1 & 0 \end{smallmatrix}\right)$ and $x_6=\left( \begin{smallmatrix} 1 & -1 \\ 1 & 0
\end{smallmatrix}\right)$ (see Example I.5.2 in~\cite{DD}). Since one can algorithmically go from
matrices to words in the presentation, and viceversa, both models are algorithmically equivalent (and we shall use one
or the other whenever convenient).

\begin{cor}\label{cosets}
The problem CIP($GL_2(\mathbb{Z})$) is solvable.
\end{cor}

\demo  There is a natural epimorphism $\varphi \colon GL_2(\mathbb{Z})\twoheadrightarrow D_{12}=\langle t_{12},\,
x_{12} \mid t_{12}^2,\, x_{12}^{12},\, (t_{12}x_{12})^2 \rangle$ given by $t_4=t_6 \mapsto t_{12}$, $x_4 \mapsto
x_{12}^3$, $x_6 \mapsto x_{12}^2$. With a few straightforward calculations, one can show that $K=\ker \varphi$ is free
of rank 2, with basis $P=[x_6,\, x_4]=\left( \begin{smallmatrix} 1 & 1 \\ 1 & 2
\end{smallmatrix}\right)$ and $Q=[x_6^2,\, x_4]=\left(
\begin{smallmatrix} 2 & 1 \\ 1 & 1 \end{smallmatrix}\right)$. Since
$K$ has index $|D_{12}|=24$ in $L=GL_2(\mathbb{Z})$, Propositions~\ref{cip-free} and~\ref{cip-fi} conclude the
proof.\qed

Finally we see that, in general, orbit decidability goes up to finite index.

\begin{prop}\label{fi}
Let $F$ be a group given in an algorithmic way, and let $A\leqslant B\leqslant Aut(F)$ be two subgroups given by finite
sets of generators, such that $A$ has finite index in $B$ and MP($A,B$) is solvable. If $A\leqslant Aut(F)$ is orbit
decidable, then $B\leqslant Aut(F)$ is orbit decidable.
\end{prop}

\demo With a coset enumeration argument similar to that in the proof of Proposition~\ref{cip-fi}, we can compute a
(finite) list, say $\{ \beta_1,\ldots ,\beta_m \}$, of left coset representatives of $A$ in $B$, i.e. $B=\beta_1
A\sqcup \cdots \sqcup \beta_m A$ (we formally need MP($A,B$) here, because $B$ may have infinite index in $Aut(F)$).
Now, for any given $u,v\in F$, the existence of $\beta \in B$ such that $u\beta$ is conjugate to $v$ is equivalent to
the existence of $i=1,\ldots ,m$ and $\alpha \in A$ such that $(u\beta_i )\alpha$ is conjugate to $v$. Hence, orbit
decidability for $B\leqslant Aut(F)$ follows from orbit decidability for $A\leqslant Aut(F)$. \qed

\subsection{The free abelian case}\label{free-abelian}

Let us concentrate on those short exact sequences in Theorem~\ref{main} with $F$ being free abelian, say
$F=\mathbb{Z}^n$, and look for orbit decidable subgroups of $Aut(\mathbb{Z}^n)=GL_n(\mathbb{Z})$.

To begin, it is an elementary fact in linear algebra that two vectors $u,v\in \mathbb{Z}^n$ lie in the same
$GL_n(\mathbb{Z})$-orbit if and only if the highest common divisor of their entries coincide and, in this case, with
the help of Euclid's algorithm, one can find an invertible matrix $A$ such that $uA=v$. In other words,

\begin{prop}\label{folk}
The full automorphism group $GL_n(\mathbb{Z})$ of a finitely generated free abelian group $\mathbb{Z}^n$, is orbit
decidable. \qed
\end{prop}

\begin{cor}
Let $A_1,\ldots ,A_m \in GL_n(\mathbb{Z})$. If $\langle A_1,\ldots ,A_m\rangle =GL_n(\mathbb{Z})$, then the
$\mathbb{Z}^n$-by-$F_m$ group $G=\mathbb{Z}^n \rtimes_{A_1,\ldots ,A_m} F_m$ has solvable conjugacy problem. \qed
\end{cor}

It is also a straightforward exercise in linear algebra to see that cyclic subgroups of $GL_n(\mathbb{Z})$ are orbit
decidable. That is, given $A\in GL_n(\mathbb{Z})$ and $u,v\in \mathbb{Z}^n$ one can algorithmically decide whether
$uA^k =v$ for some integer $k$: in fact, think $A$ as a complex matrix, work out its Jordan form (approximating
eigenvalues with enough accuracy) and then solve explicit equations (with the appropriate accuracy). This provides a
solution to the conjugacy problem for cyclic extensions of $\mathbb{Z}^n$.

\begin{prop}\label{c-ab}
Cyclic subgroups of $GL_n(\mathbb{Z})$ are orbit decidable. \qed
\end{prop}

\begin{cor}
$\mathbb{Z}^n$-by-$\mathbb{Z}$ groups have solvable conjugacy problem. \qed
\end{cor}

However, this was already known via an old result due to V.N. Remeslennikov, because $\mathbb{Z}^n$-by-$\mathbb{Z}$
groups are clearly polycyclic. In~\cite{Rem} it was proven that polycyclic groups $G$ are conjugacy separable. As a
consequence, such a group, when given by an arbitrary finite presentation, has solvable conjugacy problem (use a brute
force algorithm for solving CP$^+$($G$), and another for CP$^-$($G$) enumerating all maps into finite symmetric groups
(i.e. onto finite groups) and checking whether the images of the given elements are conjugated down there). But,
furthermore, this result can now be used to prove a more general fact about orbit decidability in $GL_n(\mathbb{Z})$.

J. Tits~\cite{Tits} proved the deep and remarkable fact that every finitely generated subgroup of $GL_n(\mathbb{Z})$ is
either virtually solvable or it contains a non-abelian free subgroup. The following proposition says that the first
kind of subgroups are always orbit decidable, so forcing orbit undecidable subgroups of $GL_n(\mathbb{Z})$ to contain
non-abelian free subgroups. This somehow means that orbit undecidability in $GL_n(\mathbb{Z})$ is intrinsically linked
to free-like structures.

%ONLY ABELIAN
%
%\begin{prop}\label{ab-ab}
%Abelian subgroups of $GL_n(\mathbb{Z})$ are orbit decidable.
%\end{prop}
%
%\demo Let $B$ be an abelian subgroup of $GL_n(\mathbb{Z})$ generated, say, by a given finite set of matrices,
%$A_1,\ldots ,A_m$. Consider the group $G=\mathbb{Z}^n\rtimes_{A_1,\ldots ,A_m}\mathbb{Z}^m$. Since it is clearly
%polycyclic, Remeslennikov's result~\cite{Rem} tells us that $G$ (with its natural presentation) has solvable conjugacy
%problem. Finally, by Theorem~\ref{main} $(a)\Rightarrow (c)$, we deduce that $B=\langle A_1,\ldots ,A_m\rangle
%\leqslant GL_n(\mathbb{Z})$ is orbit decidable (note that hypothesis (iii) of Theorem~\ref{main} is not satisfied in
%this case, but it is not used in the proof of implication $(a)\Rightarrow (c)$). \qed
%
%\begin{cor}\label{ab-ab-cor}
%Let $A_1,\ldots ,A_m \in GL_n(\mathbb{Z})$. If the $A_i$'s commute to each other then the $\mathbb{Z}^n$-by-$F_m$ group
%$G=\mathbb{Z}^n \rtimes_{A_1,\ldots ,A_m} F_m$ has solvable conjugacy problem. \qed
%\end{cor}

\begin{prop}\label{ab-solv}
Any virtually solvable subgroup of $GL_n(\mathbb{Z})$ is orbit decidable.
\end{prop}

\demo Let $B$ be a virtually solvable subgroup of $GL_n(\mathbb{Z})$, given by a finite generating set of matrices
$A_1,\ldots ,A_m$. By a Theorem of A.I. Mal'cev (see~\cite{Malc}, or Chapter 2 in~\cite{Seg}), every solvable subgroup
of $GL_n(\mathbb{Z})$ is polycyclic; so, $B=\langle A_1,\ldots ,A_m \rangle$ has a finite index subgroup $C\leqslant B$
which is polycyclic and, in particular, finitely presented.

Recurrently perform the following two lists: on one hand keep enumerating \emph{all} finite presentations of \emph{all}
polycyclic groups (use a similar strategy as that in the proof of Theorem~\ref{groups-tcp} above, enumerating first all
canonical presentations of such groups, and diagonally applying all possible Tietze transformations to all of them). On
the other hand, enumerate a set of pairs $(\mathcal{C}, \mathcal{M})$, where $\mathcal{C}$ is a finite set of
generators for a finite index subgroup $C$ of $B$, and $\mathcal{M}=\{ M_1, \ldots ,M_r \}$ is a finite set of matrices
such that $B=M_1\cdot C\cup \cdots \cup M_r\cdot C$, and in such a way that $C$ eventually visits \emph{all} finite
index subgroups of $B$; we can do this in a similar way as in the proof of Lemma~\ref{char-fi}: enumerate all saturated
and folded Stallings graphs with increasingly many vertices over the alphabet $\{ A_1,\ldots ,A_m \}$, and map the
corresponding finite index subgroup and finite set of coset representatives (one for each vertex) down to $B$, where
possible repetitions may happen (see~\cite{St}).

These two lists are infinite so the started processes will never end; but, while running, keep choosing an element in
each list in all possible ways, say $C'=\langle t_1, \ldots ,t_p \mid R_1,\ldots ,R_q\rangle$ and $(\mathcal{C},
\mathcal{M})$, and check whether there is an onto map from $\{t_1, \ldots ,t_p\}$ to $\mathcal{C}$ that extend to an
(epi)morphism $C'\to C$ (this just involves few matrix calculations in $GL_n(\mathbb{Z})$). Stop all the computations
when finding such a map (which we are sure it exists because some finite index subgroup $C\leqslant B$ is polycyclic,
and so isomorphic to one of the presentations in the first list). When this procedure terminates, we have got a finite
presentation $\langle t_1, \ldots ,t_p \mid R_1,\ldots ,R_q\rangle$ of a polycyclic group $C'$ and a map $C'\to C$ onto
a finite index subgroup $C\leqslant B=\langle A_1,\ldots ,A_m \rangle \leqslant GL_n(\mathbb{Z})$, for which we also
know a finite set of coset representatives $\mathcal{M}$, with possible repetitions.

Now, write down the natural presentation of the group $G=\mathbb{Z}^n\rtimes_C C'$. Since it is clearly polycyclic,
Remeslennikov's result~\cite{Rem} tells us that $G$ has solvable conjugacy problem (for instance, from the computed
presentation). Thus, by Theorem 3.1 (a) $\Rightarrow$ (c), the corresponding group of actions, $C\leqslant
GL_n(\mathbb{Z})$, is orbit decidable (note that hypothesis (iii) of Theorem 3.1 may not be satisfied in this case, but
it is not used in the proof of implication (a) $\Rightarrow$ (c)). Finally, $B\leqslant GL_n(\mathbb{Z})$ is orbit
decidable as well: given two vectors $u,v\in \mathbb{Z}^n$, deciding whether $v=uP$ for some $P\in B$ is the same as
deciding whether $v=uM_iQ$ for some $i=1,\ldots ,r$ and $Q\in C$, which reduces to finitely many claims to the orbit
decidability of $C\leqslant GL_n(\mathbb{Z})$. \qed

\begin{cor}\label{ab-solv-cor}
Let $A_1,\ldots ,A_m \in GL_n(\mathbb{Z})$. If $B=\langle A_1,\ldots ,A_m\rangle \leqslant GL_n(\mathbb{Z})$ is
virtually solvable, then the $\mathbb{Z}^n$-by-$F_m$ group $G=\mathbb{Z}^n \rtimes_{A_1,\ldots ,A_m} F_m$ has solvable
conjugacy problem. \qed
\end{cor}

%Note that, using Proposition~\ref{fi}, one can easily extend Propositions~\ref{c-ab} and~\ref{ab-ab} to virtually
%cyclic and virtually abelian subgroups of $GL_n(\mathbb{Z})$ (and so, get the corresponding corollaries).

Let us consider now finite index subgroups of $GL_n(\mathbb{Z})$.

\begin{prop}\label{fi-abelian}
Any finite index subgroup of $GL_n(\mathbb{Z})$ (given by generators) is orbit decidable.
\end{prop}

\demo Let $B\leqslant GL_n(\mathbb{Z})$ be a finite index subgroup generated by some given matrices. Take your favorite
presentation for $GL_n(\mathbb{Z})$ (see, for example, Section 3.5 of~\cite{MKS}) and write them in terms of it. With a
similar argument as in the proof of Lemma~\ref{char-fi}, we can compute generators for the subgroup
 $$
A=\bigcap_{P\in GL_n(\mathbb{Z})} (P^{-1}\cdot B\cdot P) =B\cap B^{P_1}\cap \cdots \cap B^{P_m} \unlhd
GL_n(\mathbb{Z}),
 $$
where $Id=P_0, P_1,\ldots ,P_m$ is a set of right coset representatives for $B$ in $GL_n(\mathbb{Z})$ (computable by
Todd-Coxeter algorithm, see Theorem~\ref{TC}). By Proposition~\ref{fi}, we are reduced to see that $A\unlhd
GL_n(\mathbb{Z})$ is orbit decidable.

Given $u,v\in \mathbb{Z}^n$, we have to decide whether some matrix of $A$ sends $u$ to $v$. Clearly, we can assume
$u,v\neq 0$ and check whether there exists $M\in GL_n(\mathbb{Z})$ such that $uM=v$ (see Proposition~\ref{folk}). Once
we have such $M$, the set of all those matrices carrying $u$ to $v$ is precisely $M\cdot Stab(v)$. And it is
straightforward to compute a finite generating set for the stabilizer of $v$ (it is conjugate to $Stab(1,0,\ldots
,0)=\{ \left( \begin{smallmatrix} 1 & 0 & \cdots & 0 \\ * & * & * & * \end{smallmatrix}\right) \}$). It remains to
algorithmically decide whether the intersection $A\cap (M\cdot Stab(v))$ is or is not empty; or, equivalently, whether
$M\in A\cdot Stab(v)$ holds or not. This is decidable because $A\cdot Stab(v)$ is a finite index subgroup of
$GL_n(\mathbb{Z})$ (here is where normality of $A$ is needed) with a computable set of generators; hence, MP($A\cdot
Stab(v),\,GL_n(\mathbb{Z})$) is solvable, again by Todd-Coxeter algorithm. \qed

\begin{cor}
Let $A_1,\ldots ,A_m \in GL_n(\mathbb{Z})$. If $\langle A_1,\ldots ,A_m\rangle$ has finite index in $GL_n(\mathbb{Z})$
then the $\mathbb{Z}^n$-by-$F_m$ group $\mathbb{Z}^n \rtimes_{A_1,\ldots ,A_m} F_m$ has solvable conjugacy problem.
\qed
\end{cor}

\medskip

Finally, let us concentrate on rank two.

\begin{prop}\label{gl2}
Every finitely generated subgroup of $GL_2(\mathbb{Z})$ is orbit decidable.
\end{prop}

\demo Let $A_1,\ldots ,A_r \in GL_n(\mathbb{Z})$ be some given matrices and consider the subgroup they generate,
$\langle A_1,\ldots ,A_r \rangle \leqslant GL_n(\mathbb{Z})$. For $n=2$, given $u,v\in \mathbb{Z}^n$, let us decide
whether there exists $A\in \langle A_1,\ldots ,A_r \rangle$ such that $uA=v$. We can clearly assume $u,v\neq 0$.

By Proposition~\ref{folk}, we can decide whether there exists $M\in GL_n(\mathbb{Z})$ such that $uM=v$ and, in the
affirmative case, find such an $M$. And it is straightforward to find a set of generators for the stabilizer of $v$,
$Stab(v)=\{ B\in GL_n(\mathbb{Z}) \mid vB=v\}$, say $\{ B_1, \ldots ,B_s \}$ (in the case $n=2$, every such stabilizer
is conjugate to $Stab(1,0)=\langle \left(
\begin{smallmatrix} 1 & 0
\\ 1 & 1
\end{smallmatrix}\right),\left( \begin{smallmatrix} 1 & 0 \\ 0 & -1 \end{smallmatrix}\right)\rangle$).
Then, the matrices sending $u$ to $v$ are precisely those contained in the coset $M\langle B_1, \ldots ,B_s\rangle$.
So, it remains to decide whether $\langle A_1, \ldots, A_r\rangle \cap M\langle B_1, \ldots ,B_s\rangle$ is empty or
not.

In the case $n=2$, this can be done algorithmically (see Corollary~\ref{cosets}). \qed

Note that the proof of Proposition~\ref{gl2} works for every dimension $n$ except at the end, when
Corollary~\ref{cosets} is used. We shall refer to this fact later.

\begin{cor}
All $\mathbb{Z}^2$-by-[f.g.$\,$free] groups have solvable conjugacy problem.\qed
\end{cor}

\subsection{The free case}\label{free}

Following the same route as in the previous subsection, let us concentrate now on those short exact sequences in
Theorem~\ref{main} with $F$ being free, say $F=F_n$, and look for orbit decidable subgroups of $Aut(F_n)$.

To begin, classical Whitehead algorithm (see Proposition~I.4.19 in~\cite{LS}) decides, given $u,v\in F_n$, whether
there exists an automorphism of $F_n$ sending $u$ to $v$ up to conjugacy. In other words,

\begin{thm}[Whitehead, \cite{Wh}]\label{Whitehead}
The full automorphism group $Aut(F_n )$ of a finitely generated free group $F_n$, is orbit decidable. \qed
\end{thm}

\begin{cor}
Let $F_n$ be a finitely generated free group. If $\varphi_1,\ldots ,\varphi_m$ generate $Aut(F_n)$, then the
$F_n$-by-$F_m$ group $G=F_n \rtimes_{\varphi_1,\ldots ,\varphi_m} F_m$ has solvable conjugacy problem. \qed
\end{cor}

Like in the abelian case, cyclic subgroups of $Aut(F_n)$ are orbit decidable by a result of P. Brinkmann. This is the
analog of Proposition~\ref{c-ab} for the free case, but here the proof is much more complicated, making strong use of
the theory of train-tracks. This was already used to solve the conjugacy problem for free-by-cyclic groups:

\begin{thm}[Brinkmann, \cite{Br}]\label{brinkmann}
Cyclic subgroups of $Aut(F_n)$ are orbit decidable. \qed
\end{thm}

\begin{cor}[Bogopolski-Martino-Maslakova-Ventura, \cite{BMMV}]
[f.g.$\,$free]-by-cyclic groups ha\-ve solvable conjugacy problem. \qed
\end{cor}

The analog of Proposition~\ref{ab-solv} and Corollary~\ref{ab-solv-cor} in the free setting is not known, and seems to
be an interesting and much more complicated problem. See Question~5 in the last section for some comments about it, and
a clear relation with Tits alternative for $Aut(F_n)$.

Let us now consider finite index subgroups of $Aut(F_n)$.

\begin{prop}\label{fi-free}
Let $F_n$ be a finitely generated free group. Any finite index subgroup of $Aut(F_n )$ (given by generators) is orbit
decidable.
\end{prop}

\demo Let $B\leqslant Aut(F_n)$ be a finite index subgroup generated by some given automorphisms. Consider Nielsen's
presentation of $Aut(F_n)$ (see Proposition~N1 in Section~3.5 of~\cite{MKS}) and write them in terms of this
presentation (i.e. as products of Nielsen automorphisms). Then, with a similar argument as in the proof of
Lemma~\ref{char-fi}, we can compute generators for the subgroup
 $$
A=\bigcap_{\phi\in Aut(F_n )} (\phi^{-1}B\phi) =B\cap B^{\phi_1}\cap \cdots \cap B^{\phi_m} \unlhd Aut(F_n),
 $$
where $Id=\phi_0, \phi_1,\ldots ,\phi_m$ is a set of right coset representatives for $B$ in $Aut(F_n)$ (computable by
Todd-Coxeter algorithm, see Theorem~\ref{TC}). By Proposition~\ref{fi}, we are reduced to see that $A\unlhd Aut(F_n)$
is orbit decidable.

Let $u,v\in F_n$. Using Whitehead's algorithm, we can check whether there exists an automorphism $\alpha\in Aut(F_n)$
carrying $u$ to $v$. Once we have such $\alpha$, the set of all those automorphisms carrying $u$ to a conjugate of $v$
is precisely $\alpha \cdot Stab(v) \cdot Inn(F_n)$. By McCool's algorithm (see Proposition~I.5.7 in~\cite{LS}), we can
compute a finite generating set for the stabilizer of $v$. It remains to algorithmically decide whether the
intersection $A\cap (\alpha \cdot Stab(v)\cdot Inn(F_n))$ is or is not empty; or, equivalently, whether $\alpha\in
A\cdot Stab(v)\cdot Inn(F_n)$ holds or not. This is decidable because $A\cdot Stab(v)\cdot Inn(F_n)$ is a finite index
subgroup of $Aut(F_n)$ (here is where normality of $A$ is needed) with a computable set of generators; hence,
MP($A\cdot Stab(v)\cdot Inn(F_n),\,Aut(F_n)$) is solvable, again by Todd-Coxeter algorithm. \qed

\begin{cor}
Let $F_n$ be a finitely generated free group. If $\varphi_1,\ldots ,\varphi_m$ generate a finite index subgroup of
$Aut(F_n)$, then the $F_n$-by-$F_m$ group $G=F_n \rtimes_{\varphi_1,\ldots ,\varphi_m} F_m$ has solvable conjugacy
problem. \qed
\end{cor}

\medskip

Now, let us concentrate on rank two. Like in the abelian case, we have

\begin{prop}\label{f2byfree}
Let $F_2$ be the free group of rank two. Then every finitely generated subgroup of $Aut(F_2)$ is orbit decidable.
\end{prop}

\demo Let $A$ be a finitely generated subgroup of $Aut(F_n)$ and let $u,v\in F_n$ be given. For $n=2$, we have to
decide whether there exists $\varphi\in A$ such that $u\varphi$ is conjugate to $v$.

Mimicking the proof of Proposition~\ref{gl2}, let us apply Whitehead's algorithm to $u,v$ (see Proposition~I.4.19
in~\cite{LS}). If there is no automorphism in $Aut(F_n)$ sending $u$ to $v$ then, clearly, the answer to our problem is
also negative. Otherwise, we have found $\alpha\in Aut(F_n)$ such that $u\alpha =v$. Now, the set of all such
automorphisms of $F_n$ is $\alpha \cdot Stab(v)$. And the set of all automorphisms of $F_n$ mapping $u$ to a conjugate
of $v$ is $\alpha \cdot Stab(v)\cdot Inn(F_n)$. By McCool's algorithm (see~ Proposition~I.5.7 in~\cite{LS}), we can
find a finite system of generators for $Stab(v)\leqslant Aut(F_n)$. Finally, all we need is to verify whether $A\cap
(\alpha\cdot Stab(v)\cdot Inn(F_n ))$ is or is not empty.

In the case $n=2$ this can be done algorithmically: since the kernel of the canonical projection
$\overline{\phantom{a}}\colon Aut(F_2 ) \twoheadrightarrow GL_2(\mathbb{Z})$ is $Inn(F_2)$ (which is a very special
fact of the rank 2 case), our goal is equivalent to verifying whether $\overline{A}\cap (\overline{\alpha}\cdot
\overline{Stab(v)})$ is or is not empty in $GL_2(\mathbb{Z})$. This can be done by Corollary~\ref{cosets}. \qed

Note that the proof of Proposition~\ref{f2byfree} works for every rank $n$ except for the last paragraph, exactly like
in Proposition~\ref{gl2}. We shall refer to this fact later.

\begin{cor}
All $F_2$-by-[f.g.$\,$free] groups have solvable conjugacy problem. \qed
\end{cor}

Another nice examples of orbit decidable subgroups in $Aut(F_n)$ come from geometry. Certain mapping class groups of
surfaces with boundary and punctures turn out to embed in the automorphism group of the free group of the appropriate
rank. The image of these embeddings are easily seen to be orbit decidable in two special cases. From~\cite{DF} we
extract the following two examples.

Let $S_{g,b,n}$ be an orientable surface of genus $g$, with $b$ boundary components and $n$ punctures. It is well known
that its fundamental group has presentation
 $$
\Sigma_{g,b,n}=\langle x_1,y_1,\ldots ,x_g,y_g,\, z_1,\ldots ,z_b,\, t_1,\ldots ,t_n \mid [x_1,y_1]\cdots
[x_g,y_g]z_1\cdots z_b t_1\cdots t_n =1\rangle ,
 $$
and, except for $b=n=0$, is a free group of rank $2g+b+n-1$. In the following cases, the mapping class group of
$S_{g,b,n}$ can be viewed as a subgroup of $Aut(F_{2g+b+n-1})$, see~\cite{DF} for details:

\begin{prop}(see~\cite{DF})
Let $S_{g,b,n}$ be an orientable surface of genus $g$, with $b$ boundary components, and with $n$ punctures.
\begin{itemize}
\item[(i)] {\rm (\textbf{Maclachlan})} The positive pure mapping class group of $S_{g,0,n+2}$ becomes (when the basepoint
is taken to be the ($n+2$)nd puncture) the group $Aut^+_{g,0,1^{\perp n+1}}$ of automorphisms of
$\Sigma_{g,0,n+1}\simeq F_{2g+n}$ which fix each conjugacy class $[t_j]$, $j=1,\ldots ,n+1$ (the case $g=0$ gives the
pure braid group on $n+1$ strings modulo the center, $B_{n+1}/Z(B_{n+1})$).

\item[(ii)] {\rm (\textbf{A'Campo})} The positive mapping class group of $S_{g,1,n}$ becomes (when the basepoint is taken
to be on the boundary) the group $Aut^+_{g, \hat 1, n}$ of automorphisms of $\Sigma_{g,1,n}\simeq F_{2g+n}$ which fix
$z$ and permute the set of conjugacy classes $\{ [t_1],\ldots ,[t_n]\}$ (the case $g=0$ gives the braid group on $n$
strings, $B_n$).
\end{itemize}
\end{prop}

A particularly interesting case is when $g=0$ in (ii) above: $Aut^+_{0, \hat 1, n}$ is the image of classical Artin's
embedding of the braid group on $n$ strings into $Aut(F_n)$, sending generator $\sigma_i \in B_n$ ($i=1,\ldots ,n-1$)
to $\sigma_i \colon F_n \to F_n$, $t_i\mapsto t_it_{i+1}t_i^{-1}$, $t_{i+1}\mapsto t_i$, $t_j\mapsto t_j$ for $j\neq
i,i+1$. The subgroup $Aut^+_{0, \hat 1, n} =\langle \sigma_1,\ldots ,\sigma_{n-1}\rangle \leqslant Aut(F_n)$ is then
characterized as those automorphisms $\varphi \in Aut(F_n)$ for which $(t_1t_2\cdots t_n)\varphi =t_1t_2\cdots t_n$ and
there exist words $w_1,\ldots ,w_n\in F_n$ and a permutation $\sigma$ of the set of indices such that  $t_i\varphi =
w_i^{-1}t_{\sigma(i)}w_i$.

All these groups of automorphisms, $Aut^+_{g,0,1^{\perp n+1}}\leqslant Aut(F_{2g+n})$ and $Aut^+_{g, \hat 1,
n}\leqslant Aut(F_{2g+n})$, are easily seen to be orbit decidable because of the following observation.

\begin{prop}\label{odmcg}
Let $F_n$ be a finitely generated free group, and let $u_i,v_i \in F_n$ ($i=0,\ldots ,m$) be two lists of elements.
Then,
\begin{itemize}
\item[(i)] $A=\{ \varphi \mid u_0\varphi=v_0,\, u_1\varphi \sim v_1,\ldots ,u_m\varphi \sim v_m \}\leqslant Aut(F_n)$ and
\item[(ii)] $B=\{ \varphi \mid u_0\varphi=v_0,\, u_1\varphi \sim v_{\sigma(1)},\ldots ,u_m\varphi \sim v_{\sigma(m)}
\text{ for some }\sigma\in Sym(m)\}\leqslant Aut(F_n)$
\end{itemize}
are orbit decidable.
\end{prop}

\demo For (i), given $u,v\in F_n$, we have to decide whether there exists an automorphism $\varphi\in Aut(F_n)$ such
that $u_0\varphi =v_0$, $[u_1]\varphi =[v_1],\ldots ,[u_m]\varphi =[v_m]$ and $[u]\varphi =[v]$, where brackets denote
conjugacy classes. This is the same as deciding whether there exists $\varphi$ such that $[u_0]\varphi =[v_0]$,
$[u_1]\varphi =[v_1],\ldots ,[u_m]\varphi =[v_m]$ and $[u]\varphi =[v]$ (and, in the affirmative case, composing
$\varphi$ by $\gamma_{w^{-1}_0}$, where $w_0$ is the first conjugator, $u_0\varphi =w_0^{-1}v_0w_0$). One can make this
decision by applying Proposition~4.21 in~\cite{LS}. Finally, (ii) can be solved by using up to $m!$ many times the
solution given for (i). \qed

It is worth mentioning that D. Larue in his PhD thesis~\cite{Larue} analyzed the $Aut^+_{0, \hat 1, n}$-orbit of $t_1$
in $\Sigma_{0,1,n}$ (i.e. the $B_n$-orbit of $t_1$ in $F_n$) and he provided an algorithm to decide whether a given
word $w\in F_n$ belongs to this orbit (note that this is not exactly a special case of OD($B_n$) because $B_n$ does not
contain all inner automorphisms). Although working only for the orbit of $t_1$, the algorithm provided is faster and
nicer than that provided in Proposition~\ref{odmcg}.

\begin{cor}
Let $F_{2g+n}$ be a finitely generated free group. If $\varphi_1,\ldots ,\varphi_m \in Aut(F_{2g+n})$ generate the
positive pure mapping class group $Aut^+_{g,0,1^{\perp n+1}}$, or the positive mapping class group $Aut^+_{g, \hat 1,
n}$ then the $F_{2g+n}$-by-$F_m$ group $G=F_{2g+n} \rtimes_{\varphi_1,\ldots ,\varphi_m} F_m$ has solvable conjugacy
problem (a particular case of this being when $\varphi_1,\ldots ,\varphi_m$ generate the standard copy of the braid
group $B_n\leqslant Aut(F_n)$). \qed
\end{cor}

\section{Negative results}\label{-}

Let us construct now negative examples, namely orbit undecidable subgroups of $GL_n(\mathbb{Z})$ and $Aut(F_n)$ which,
of course, will correspond to [f.g.$\,$free abelian]-by-[f.g.$\,$free] and [f.g.$\,$free]-by-[f.g.$\,$free] groups with
unsolvable conjugacy problem.

As mentioned above, Miller constructed a [f.g.$\,$free]-by-[f.g.$\,$free] group with unsolvable conjugacy problem
(see~\cite{M1}); here, we have a first source of examples of finitely generated subgroups of the automorphism group of
a free group, which are orbit undecidable. In the present section, we will generalize this construction by giving a
source of orbit undecidability in $Aut(F)$ for more groups $F$. When taking $F$ to be free, this will reproduce
Miller's example; when taking $F=\mathbb{Z}^n$ for $n\geqslant 4$, we will obtain orbit undecidable subgroups in
$GL_n(\mathbb{Z})$, which correspond to the first known examples of [f.g.$\,$free abelian]-by-[f.g.$\,$free] groups
with unsolvable conjugacy problem.

Let us recall Miller's construction. It begins with an arbitrary finite presentation, $H=\langle s_1,\dots ,s_n \mid
R_1,\dots ,R_m\rangle$, where the $R_j$'s are words on the $s_i$'s. Let $F_{n+1}=\langle q,s_1,\dots ,s_n \mid \,\,
\rangle$ and $F_{m+n}=\langle t_1,\dots ,t_m,d_1,\dots ,d_n \mid \,\, \rangle$ be the free groups of rank $n+1$ and
$m+n$, respectively, on the listed generators. Consider now the $m+n$ automorphisms of $F_{n+1}$ given by
 $$
\begin{array}{rclcrcl}
\alpha_i \colon F_{n+1} & \to & F_{n+1} & & \beta_j \colon F_{n+1} & \to & F_{n+1} \\ q & \mapsto & qR_i & \quad \quad
& q & \mapsto & s_j^{-1}qs_j \\ s_k & \mapsto & s_k & & s_k & \mapsto & s_k
\end{array} ,
 $$
for $i=1,\ldots ,m$ and $j,k=1,\ldots ,n$, and denote the group of automorphisms they generate by $A(H)\leqslant
Aut(F_{n+1})$. Next, consider the $F_{n+1}$-by-$F_{m+n}$ group defined by these automorphisms,
 $$
G(H)=F_{n+1}\rtimes_{\alpha_1,\ldots ,\alpha_m, \beta_1,\ldots ,\beta_n} F_{m+n}.
 $$
% t_j^{-1}qt_j & = & qR_j, & \quad \quad &
%d_k^{-1}qd_k & = & s_k^{-1}qs_k, \\ t_j^{-1}s_i t_j & = & s_i, & & d_k^{-1}s_id_k & = & s_i,
The following Theorem is Corollary~5 in Chapter~III of~\cite{M1}. Below, we shall provide an alternative proof.

\begin{thm}[Miller, \cite{M1}]\label{Miller}
If $H$ has unsolvable word problem then $G(H)$ has unsolvable conjugacy problem.
\end{thm}

So, applying Miller's construction to a presentation $H$ with $n$ generators, $m$ relations, and with unsolvable word
problem, one obtains a $(m+n)$-generated subgroup of $Aut(F_{n+1})$, namely $A(H)$, which is orbit undecidable.

In \cite{Bor}, V. Borisov constructed a group presented with 4 generators, 12 relations, and having unsolvable word
problem. In order to reduce the number of generators to 2 (and so have the corresponding orbit undecidable subgroup
living inside $Aut(F_3)$) we can use Higman-Neumann-Neumann embedding theorem, saying that any countable group can be
embedded in a group with two generators and the same number of relations (see~\cite{HNN}). Since solvability of the
word problem clearly passes to subgroups, we obtain a group with $n=2$ generators, $m=12$ relations, and having
unsolvable word problem. Using Miller's construction we conclude the existence of a $F_3$-by-$F_{14}$ group with
unsolvable conjugacy problem. In other words,

\begin{cor}\label{14}
There exists a 14-generated subgroup $A\leqslant Aut(F_3)$ which is orbit undecidable. \qed
\end{cor}

\medskip

Let us now find a more general source of orbit undecidability that will apply to more groups $F$ other than free (and,
in the free case, will coincide with Miller's example via Mihailova's result).

Let $F$ be a group. Recall that the \emph{stabilizer} of a given subgroup $K\leqslant F$, denoted $Stab(K)$, is
 $$
Stab\,(K)=\{ \varphi \in Aut(F) \mid k\varphi =k\quad \forall k\in K\} \leqslant Aut(F).
 $$
For simplicity, we shall write $Stab\,(k)$ to denote $Stab\,(\langle k \rangle)$, $k\in F$. Furthermore, we define the
\emph{conjugacy stabilizer} of $K$, denoted $Stab^*(K)$, to be the set of automorphisms acting as conjugation on $K$,
formally $Stab^*(K)=Stab\,(K)\cdot Inn(F) \leqslant Aut(F)$.

\begin{prop}\label{key}
Let $F$ be a group. Suppose we are given two subgroups $A\leqslant B\leqslant Aut(F)$ and an element $v\in F$ such that
$B\cap Stab^*(v) =\{ Id\}$. If $A\leqslant Aut(F)$ is orbit decidable then MP($A,B$) is solvable.
\end{prop}

\demo Given $\psi \in B\leqslant Aut(F)$, let us decide whether $\psi\in A$ or not. Take $w=v\psi$ and observe that
 $$
\{ \phi \in B \mid v\phi \sim w\} =B\cap (Stab^*(v) \cdot \psi )=(B\cap Stab^*(v) )\cdot \psi =\{ \psi \}.
 $$
So, there exists $\phi \in A$ such that $v\phi$ is conjugate to $w$ in $F$, if and only if $\psi \in A$. Hence, orbit
decidability for $A\leqslant Aut(F)$ solves MP($A,B$). \qed

One can interpret Proposition~\ref{key} by saying that if, for a group $F$, $Aut(F)$ contains a pair of subgroups
$A\leqslant B\leqslant Aut(F)$ with unsolvable MP($A,B$) then $A\leqslant Aut(F)$ is orbit undecidable.

The most classical example of unsolvability of the membership problem goes back to fifty years ago. In~\cite{Mih} (see
also Chapter III.C of~\cite{M1}) Mihailova gave a nice example of unsolvability of the membership problem. The
construction goes as follows.

Like before, start with an arbitrary finite presentation, $H=\langle s_1,\dots ,s_n \mid R_1,\dots ,R_m\rangle$, and
consider the subgroup $A=\{ (x,y)\in F_n \times F_n \mid x=_{H}y\} \leqslant F_n \times F_n$. It is straightforward to
verify that $A=\langle (1,R_1),\ldots ,(1,R_m),\, (s_1,s_1),\ldots ,(s_n,s_n)\rangle$ (and so it is finitely
generated), and that MP($A,F_n \times F_n$) is solvable if and only if WP($H$) is solvable.

By Higman-Neumann-Neumann embedding theorem, we can restrict our attention to 2-gene\-ra\-ted groups (take $n=2$ in the
above paragraph). From all this, we deduce the following.

\begin{thm}\label{f2xf2}
Let $F$ be a finitely generated group such that $F_2 \times F_2$ embeds in $Aut(F)$ in such a way that the image
intersects trivially with $Stab^*(v)$, for some $v\in F$. Then, $Aut(F)$ contains an orbit undecidable subgroup; in
other words, there exist $F$-by-[f.g.$\,$free] groups with unsolvable conjugacy problem. \qed
\end{thm}

In the rest of the section, we shall use Theorem~\ref{f2xf2} to obtain explicit examples in the free abelian and free
cases.

\subsection{The free abelian case}\label{-free-abelian}

It is well known that $F_2$ embeds in $GL_2(\mathbb{Z})$ and so, $F_2\times F_2$ embeds in $GL_4(\mathbb{Z})$. Hence,
we can deduce the following result.

\begin{prop}\label{gl4}
For $n\geqslant 4$, $GL_n(\mathbb{Z})$ contains finitely generated orbit undecidable subgroups.
\end{prop}

\demo Consider the subgroup of $GL_2(\mathbb{Z})$ generated by $P=\left( \begin{smallmatrix} 1 & 1
\\ 1 & 2
\end{smallmatrix}\right)$ and $Q=\left( \begin{smallmatrix} 2 & 1 \\ 1 & 1 \end{smallmatrix}\right)$, which is free and
freely generated by $\{ P,Q\}$ as discussed in the proof of Corollary~\ref{cosets}. We claim that $\langle P,Q\rangle
\cap
Stab^*\big( (1,0)\big) =\langle \left( \begin{smallmatrix} 1 & 0 \\
12 & 1 \end{smallmatrix}\right) \rangle$. In fact, it is clear that $Stab^*\big( (1,0)\big)=Stab\,\big( (1,0)\big)
=\left\{ \left( \smallmatrix 1 & \phantom{\pm} 0 \\ n & \pm 1
\endsmallmatrix \right) \mid n\in \mathbb{Z} \right\}$ (and we can forget the negative signum
because we are interested in the intersection with $\langle P,Q\rangle \leqslant SL_2(\mathbb{Z})$). Now, the image of
$\left( \smallmatrix 1 & 0 \\ 1 & 1 \endsmallmatrix \right) =x_6^{-1}x_4$ under $\varphi$ is
$x_{12}^{-2}x_{12}^3=x_{12}\in D_{12}$ (see the proof of Corollary~\ref{cosets} for notation). So,
 $$
\langle P,Q\rangle \cap Stab^*\big( (1,0)\big) =\ker \varphi \cap Stab \big( (1,0)\big) =\langle
(x_6^{-1}x_4)^{12}\rangle =\langle \left( \begin{smallmatrix} 1 & 0 \\ 12 & 1\end{smallmatrix}\right) \rangle.
 $$
Choose now a (free) subgroup $\langle P',Q'\rangle \leqslant \langle P,Q\rangle$ intersecting trivially with the cyclic
subgroup $\langle \left( \begin{smallmatrix} 1 & 0 \\ 12 & 1\end{smallmatrix}\right) \rangle$ (this always exists in
non-cyclic free groups). And, for $n\geqslant 4$, consider
 $$
B=\left\langle \left( \begin{array}{c|c} P' & 0 \\ \hline 0 & Id
\end{array} \right) ,\left(
\begin{array}{c|c} Q' & 0 \\ \hline 0 & Id \end{array} \right) ,\left( \begin{array}{c|c} Id & 0 \\ \hline 0 & P'
\end{array} \right) ,\left( \begin{array}{c|c} Id & 0 \\ \hline 0 & Q' \end{array} \right) \right\rangle \leqslant
GL_4(\mathbb{Z}) \leqslant GL_n(\mathbb{Z}),
 $$
which is clearly isomorphic to $F_2 \times F_2$. By construction, $B$ intersects trivially with the (conjugacy)
stabilizer of $v=(1,0,1,0,\ldots ,0)\in \mathbb{Z}^n$. Finally, using Mihailova's construction, find a finitely
generated subgroup $A\leqslant B$ with unsolvable MP($A,B$). By Proposition~\ref{key} applied to $F=\mathbb{Z}^n$, $A$
is a finitely generated orbit undecidable subgroup of $Aut(\mathbb{Z}^n)=GL_n(\mathbb{Z})$. \qed

\begin{cor}
There exist $\mathbb{Z}^4$-by-[f.g.$\,$free] groups with unsolvable conjugacy problem. \qed
\end{cor}

While the constructions are quite different, these groups are reminiscent of Miller's examples, but with a free abelian
base group. To the best of our knowledge, they are the first known examples of [f.g.$\,$free abelian]-by-[f.g.$\,$free]
groups with unsolvable conjugacy problem. As a side consequence of the previous reasoning, we also obtain the following
corollary.

\begin{cor}\label{cipgl4}
For $n\geqslant 4$, CIP($GL_n(\mathbb{Z})$) is unsolvable.
\end{cor}

\demo As noted above, the proof of Proposition~\ref{gl2} works entirely for any dimension $n$ for which
CIP($GL_n(\mathbb{Z})$) is solvable (for example, $n=2$). But Proposition~\ref{gl4} states the existence of finitely
generated orbit undecidable subgroups of $GL_n(\mathbb{Z})$, for $n\geqslant 4$. Hence, CIP($GL_n(\mathbb{Z})$) must be
unsolvable in this case. \qed

\subsection{The free case}\label{-free}

In order to apply Theorem~\ref{f2xf2} to the free group $F_3 =\langle q,a,b \mid \quad\rangle$ of rank 3, we need to
identify a copy of $F_2 \times F_2$ inside $Aut(F_3)$. For every $w\in \langle a,b\rangle$, consider the automorphisms
$_w\theta_1 \colon F_3 \to F_3$, $q\mapsto wq$, $a\mapsto a$, $b\mapsto b$, and $_1 \theta_w \colon F_3 \to F_3$,
$q\mapsto qw$, $a\mapsto a$, $b\mapsto b$. Clearly, $_{w_1}\theta_1 \, _{w_2}\theta_1 =\, _{w_1w_2}\theta_1$ and $_1
\theta_{w_1} \, _1\theta_{w_2} =\, _1\theta_{w_2w_1}$, which means that $\{ \, _w\theta_1 \mid w\in \langle
a,b\rangle\} \simeq F_2$ and $\{ \, _1\theta_w \mid w\in \langle a,b\rangle\} \simeq F_2^{\rm op}\simeq F_2$. It is
also clear that $_{w_1}\theta_1 \, _1\theta_{w_2} =\, _{w_1}\theta_{w_2} =\, _1\theta_{w_2} \, _{w_1}\theta_1$ (with
the natural definition for $_{w_1}\theta_{w_2}$). So, we have an embedding $F_2\times F_2 \simeq F_2^{\rm op}\times
F_2^{\rm op} \to Aut(F_3)$ given by $(w_1, w_2)\mapsto \, _{w_1^{-1}}\theta_{w_2}$, whose image is
 $$
B=\langle _{a^{-1}}\theta_1 ,\, _{b^{-1}}\theta_1,\, _1\theta_a ,\, _1\theta_b \rangle =\{ \, _{w_1}\theta_{w_2} \mid
w_1,\, w_2\in \langle a,b\rangle\} \leqslant Aut(F_3).
 $$
In order to use Proposition~\ref{key}, let us consider the element $v=qaqbq$. We claim that $B\cap Stab^*(v) =\{ Id\}$.
In fact, suppose $w_1,w_2\in \langle a,b\rangle$ are such that $(v)\, _{w_1}\theta_{w_2}=w_1qw_2aw_1qw_2bw_1qw_2$ is
conjugate to $v=qaqbq$ in $F_3$. Since both words have exactly three occurrences of $q$, they must agree up to cyclic
reordering. That is, $q(w_2aw_1)q(w_2bw_1 )q(w_2w_1 )$ equals either $qaqbq$, or $qbq^2a$, or $q^2aqb$. From this, one
can straightforward deduce that $w_1=w_2=1$ in all three cases. Thus, $_{w_1}\theta_{w_2}=Id$ proving the claim.

Now, let $H=\langle a,b \mid R_1,\ldots ,R_{12}\rangle$ be Borisov's example of a group with unsolvable word problem,
embedded in a 2-generated group via Higman-Neumann-Neumann embedding (see above and~\cite{HNN}). By Mihailova result
and Proposition~\ref{key}, $A=\langle \, _1\theta_{R_1},\, \ldots ,\, _1\theta_{R_{12}} ,\, _{a^{-1}}\theta_{a},\,
_{b^{-1}}\theta_{b}\rangle \leqslant Aut(F_3)$ is orbit undecidable. Hence, by Theorem~\ref{freebyfree}, the
$F_3$-by-$F_{14}$ group determined by the automorphisms $\, _1\theta_{R_1},\, \ldots ,\, _1\theta_{R_{12}},$
$_{a^{-1}}\theta_a,\, _{b^{-1}}\theta_b \in Aut(F_3)$,
  $$
G=\left\langle \begin{array}{c|ccccc}  & t_i^{-1}qt_i =qR_i & & d_1^{-1}qd_1 =a^{-1}qa & & d_2^{-1}qd_2 =b^{-1}qb
\\ q,a,b,t_1,\ldots ,t_{12},d_1,d_2 & t_i^{-1}at_i =a & , & d_1^{-1}ad_1 =a & , & d_2^{-1}ad_2 =a \\ & t_i^{-1}bt_i =b & &
d_1^{-1}bd_1 =b & & d_2^{-1}bd_2 =b \end{array} \right\rangle ,
 $$
has unsolvable conjugacy problem. This is precisely Miller's group $G(H)$ associated to $H=\langle a,b \mid R_1,\ldots
,R_{12}\rangle$ (see the beginning of the present section). Thus, the argument just given provides an alternative proof
of Miller's Theorem~\ref{Miller}.

\begin{cor}\label{autf3}
For $n\geqslant 3$, CIP($Aut(F_n)$) is unsolvable.
\end{cor}

\demo  As noted above, the proof of Proposition~\ref{f2byfree} works entirely for any rank $n$ for which
CIP($Aut(F_n)$) is solvable. But, for $n\geqslant 3$, $Aut(F_n)$ contains finitely generated orbit undecidable
subgroups. Hence, CIP($Aut(F_n)$) must be unsolvable in this case. \qed

\section{Open problems}\label{open}

Finally, we collect several questions suggested by the previous results.

\bigskip

\noindent {\bf Question 1.} Apart from finitely generated abelian, free, surface and polycyclic groups, and virtually
all of them (see Theorem~\ref{groups-tcp}), find more examples of groups $F$ with solvable twisted conjugacy problem.

\medskip \noindent{\it Commentary.} As mentioned in Section~\ref{more-applic}, for every group $F$ with solvable twisted
conjugacy problem, the study of orbit decidability/undecidability among subgroups of $Aut(F)$ becomes interesting
because it directly corresponds to solving the conjugacy problem for some extensions of $F$.\qed

\bigskip

\noindent {\bf Question 2.} Is the twisted conjugacy problem solvable for finitely generated hyperbolic groups ?

\medskip \noindent{\it Commentary.} The first step into this direction is the solvability of the twisted conjugacy
problem for finitely generated free groups, proven in~\cite{BMMV}. However, there is no hope to extend that proof for
hyperbolic groups because we do not have enough control on the automorphism group of an arbitrary hyperbolic group.
\qed

\bigskip

\noindent {\bf Question 3.} Let $F$ be a group given by a finite presentation $\langle X\mid R \rangle$, and suppose we
are given a set of words $\{ w_1,\ldots, w_r \}$ on $X$ such that $K=\langle w_1,\ldots ,w_r \rangle\leqslant F$ is a
finite index subgroup. Does solvability of TCP($F$) imply solvability of TCP($K$) ? Is it true with the extra
assumption that $K$ is characteristic in $F$ ?

\medskip \noindent{\it Commentary.} The reverse implication is proved to be true in Proposition~\ref{R-fi}~(i), under the
characteristic assumption for $K$. However, to go down from $F$ to $K$ we would have to consider the apparently more
complicated problem of dealing with possible automorphism of $K$ which do not extend to automorphisms of $F$. Maybe
this is a strong enough reason to build a counterexample. Note that the answer to the non-twisted version of the same
question is \emph{no} by a result of Collins-Miller (see~\cite{CM2}). \qed

\bigskip

\noindent {\bf Question 4.} Let $F$ be $\mathbb{Z}^n$ or $F_n$, and let $A\leqslant B\leqslant Aut(F)$ be two subgroups
given by finite sets of generators, such that $A$ has finite index in $B$, and MP($A,B$) is solvable. Is it true that
orbit decidability for $B\leqslant Aut(F)$ implies orbit decidability for $A\leqslant Aut(F)$ ?

\medskip \noindent{\it Commentary.} The reverse implication is proven in Proposition~\ref{fi}. With the first argument there, a finite list of left coset
representatives of $A$ in $B$ can be computed, say $B=\beta_1 A\sqcup \cdots \sqcup \beta_m A$. Then, given $u,v\in F$,
and assuming we got $\beta\in B$ such that $v\sim u\beta$, the set of all such $\beta$'s is $\beta\cdot Stab(v)\cdot
Inn(F)$. Since generators for $Stab(v)\cdot Inn(F)$ are computable (by straightforward matrix calculations in the case
$F=\mathbb{Z}^n$, and by McCool's algorithm in the case $F=F_n$), it only remains to decide whether the intersection
$\beta\cdot Stab(v)\cdot Inn(F)\cap A$ is empty or not. This can be done in the case $n=2$ because
CIP($GL_2(\mathbb{Z})$) is solvable (see the proof of Proposition~\ref{f2byfree}). However, this last part of the
argument does not work in the other cases, because CIP($Aut(F)$) is unsolvable for $F=\mathbb{Z}^n$ with $n\geq 4$ (see
Corollary~\ref{cipgl4}), and for $F=F_n$ with $n\geq 3$ (see Corollary~\ref{autf3}).

As partial answers, note that Propositions~\ref{fi-abelian} and~\ref{fi-free} show that the answer is \emph{yes} in the
special cases where $F=\mathbb{Z}^n$ and $B=GL_n(\mathbb{Z})$, and where $F=F_n$ and $B=Aut(F_n)$, respectively.

We formulate the question for free and free abelian groups because, if one allows an arbitrary ambient $F$, then the
answer is negative: consider Collins-Miller example of a finitely presented group $G$ with an index two subgroup
$F\leqslant G$ such that CP($G$) is solvable but CP($F$) is unsolvable (see~\cite{CM2}); furthermore, $G$ contains an
element $g_0 \in G$ of order two which acts on $F$ as a non-inner automorphism $\gamma_{g_0}\in Aut(F)$. So, we have
the short exact sequence $1\to F\to G\to C_2 \to 1$. Since CP($G$) is solvable, the action subgroup $B=\{Id,
\gamma_{g_0}\}\leqslant Aut(F)$ is orbit decidable; however, CP($F$) is unsolvable meaning that the trivial subgroup
$A=\{ Id\}\leqslant Aut(F)$ is orbit undecidable. \qed

\bigskip

\noindent {\bf Question 5.} Is any virtually solvable subgroup of $Aut(F_n)$ orbit decidable?

\medskip \noindent{\it Commentary.} This is the analog of Proposition~\ref{ab-solv} in the free setting. However, it
reduces to the same question for virtually free abelian subgroups. In fact, Bestvina-Feighn-Handel proved
in~\cite{BFH1} that every solvable subgroup of $Out(F_n)$ contains a finitely generated free abelian subgroup of index
at most $3^{5n^2}$ (additionally, it is also known that free abelian subgroups of $Out(F_n)$ have rank at most $2n-3$).
And the same if true for $Aut(F_n)$ because one can easily embed $Aut(F_n)$ in $Out(F_{2n})$ by sending $\alpha\in
Aut(F_n)$ to the outer automorphism of $F_{2n}$ which acts as $\alpha$ on both the first half and the second half of
the generating set. So, the situation is formally simpler than in Proposition~\ref{ab-solv}, but the argument there
does not work here because we cannot use the trick about polycyclic groups. Apart from the possible finite index step,
this question asks for a multidimensional version of Brinkmann's result (Theorem~\ref{brinkmann}). So, due to the
complexity of the proof and solution for the cyclic case, it seems a quite difficult question.

It is worth remarking that Bestvina-Feighn-Handel also proved in~\cite{BFH2} a strong version of Tits alternative for
$Out(F_n)$: every subgroup of $Out(F_n)$ is either virtually solvable (and hence virtually free abelian) or contains a
non-abelian free group. Since the same is true for subgroups of $Aut(F_n)$ via the above embedding, an affirmative
answer to Question~5 would then force orbit undecidable subgroups of $Aut(F_n)$ to contain non-abelian free subgroups,
like in the abelian context. This would intuitively confirm that, again, orbit undecidability is intrinsically linked
to free-like structures. \qed

\bigskip

\noindent {\bf Question 6.} Is any finitely presented subgroup of $Aut(F_n)$ orbit decidable?

\medskip \noindent{\it Commentary.}  This question contains Question~5, so it is even more difficult to be answered
in the affirmative. Note that orbit undecidable subgroups of the form $A=\langle \, _1\theta_{R_1},\, \ldots ,\,
_1\theta_{R_{12}} ,\, _{a^{-1}}\theta_{a},\, _{b^{-1}}\theta_{b}\rangle \leqslant Aut(F_3)$ corresponding to Miller's
examples (see subsection~\ref{-free}) are \emph{not} a counterexample to this question because they are not finitely
presented by Proposition~B in Grunewald~\cite{Gru} (there, $F\times_{\phi}F$ corresponds to our $A\leqslant B\simeq
F_2\times F_2 \leqslant Aut(F_3)$, and $H$ corresponds to Borisov's group with two generators and unsolvable word
problem). Alternatively, $A$ is not finitely presented because it is not the direct product of finite index subgroups
of $F_2 =\{ \, _w\theta_1 \mid w\in \langle a,b\rangle\}$ and $F_2 =\{ \, _1\theta_w \mid w\in \langle a,b\rangle\}$
(see Short's description of finitely presented subgroups of $F_2\times F_2$ in~\cite{Sh}). \qed

\bigskip

\noindent {\bf Question 7.} Are there more sources of orbit undecidability other than exploiting the unsolvability of
membership problem for certain subgroups ?

\medskip \noindent{\it Commentary.} In order to find new sources, one needs to relate orbit decidability with some
other algorithmic problem, for which there are known unsolvable examples.\qed

\bigskip

\noindent {\bf Question 8.} Is it true that every finitely generated subgroup of $GL_3(\mathbb{Z})$ is orbit decidable
? Or conversely, is it true that there exists a $\mathbb{Z}^3$-by-free group with unsolvable conjugacy problem ? In
close relation with this, is CIP($GL_3(\mathbb{Z})$) solvable ?

\medskip \noindent{\it Commentary.} Propositions~\ref{gl2} and~\ref{gl4}, and Corollaries~\ref{cosets} and~\ref{cipgl4}    show that the cases of dimension 2 and dimension
bigger than or equal to 4 behave oppositely with respect to these three questions (answers being \emph{yes, no, yes},
and \emph{no, yes, no}, respectively). For the case of dimension 3, we point out that $GL_3(\mathbb{Z})$ is not
virtually free, so the argument given in Proposition~\ref{gl2} does not work in this case. But, on the other hand, $F_2
\times F_2$ does not embed in $GL_3(\mathbb{Z})$ either (in fact, only very simple groups $G$ satisfy $G\times G
\leqslant GL_3(\mathbb{Z})$), so the argument in Proposition~\ref{gl4} does not work in dimension 3 either (unless one
can find other pairs of subgroups $A\leqslant B\leqslant GL_3(\mathbb{Z})$ with unsolvable MP($A,B$)). In the free
context, the situation is easier, with the difference in behavior happening between rank two and rank three (see
Propositions~\ref{f2byfree} and~\ref{14}).

Also, it is interesting to remark that this situation is very similar (and maybe related) to the coherence of linear
groups: it is known that $GL_2(\mathbb{Z})$ is coherent, because it is virtually free, and that $GL_n(\mathbb{Z})$ is
not coherent for $n\geqslant 4$, precisely because it contains $F_2 \times F_2$ (see for example~\cite{Gru}
and~\cite{Sh}). The question is still open in dimension 3, where none of the previous arguments work.\qed

\bigskip

%\noindent {\bf Question 9.} Construct a finitely presented group $G$ with $Aut(G)$ containing an orbit undecidable
%infinite cyclic subgroup. Is it true that every finitely presented group $C$ is isomorphic to an orbit undecidable
%subgroup of $Aut(G)$ for some $G$ ?
%
%\medskip \noindent{\it Commentary.} Directly from definitions, CP($G$) is unsolvable if and only if the trivial subgroup
%$\{ Id\}\leqslant Aut(G)$ is orbit undecidable. An straightforward argument shows that if CP($G$) is unsolvable then
%the group of automorphisms permuting the coordinates of $\ast_{i=1}^p G$, say $S_p \simeq C\leqslant Aut(\ast_{i=1}^p
%G)$, (and similarly every subgroup of it) is orbit undecidable. Hence, every finite group $C$ can be realized as an
%orbit undecidable subgroup of $Aut(G)$ for appropriate $G$. The question asks for the next natural step, namely
%infinite cyclic. Or in general, it is conceivable that every finitely presented group $C$ can be realized in this way.

\noindent {\bf Question 9.} Is it true that
\begin{itemize}
\item[(a)] for every two finitely generated groups $A$ and $B$ (except for $A=1$ and $|B|<\infty$), there exists a finitely presented group $G$
such that $Aut(G)$ simultaneously contains an orbit decidable subgroup isomorphic to $A$, and an orbit undecidable
subgroup isomorphic to $B$ ?

\item[(b)] for every finitely generated group $A$, there exists a finitely presented group $G$ such that $Aut(G)$ contains
an orbit decidable subgroup isomorphic to $A$ ?

\item[(c)] for every finitely generated infinite group $B$, there exists a finitely presented group $G$ such that CP($G$) is
solvable, and $Aut(G)$ contains an orbit undecidable subgroup isomorphic to $B$ ?
\end{itemize}

\medskip \noindent{\it Commentary.} Question (b) asks for a positive decisional condition, and question (c) for the corresponding negative one, while
question (a) asks whether they are compatible within the same group $G$. Formally, (b) and (c) are partial cases of (a)
(note that CP($G$) is equivalent to OD($\{Id\}$)). However, an easy construction using the direct product shows that
affirmative answers for (b) and (c) would imply an affirmative answer for (a) too.

Before, note that $A=1$ and $|B|<\infty$ is the only situation where the copy of $A$ will necessarily be a finite index
subgroup of the copy of $B$ and so, Proposition~\ref{fi} would then say that solvability for OD($A$) implies
solvability for OD($B$). In all other cases, even if $B$ has a finite index subgroup isomorphic to $A$, it is
conceivable that $Aut(G)$ could contain copies of $A$ and $B$ apart enough to each other to fulfill the requirements of
question (a).

Now, let $G_1$ and $G_2$ be two groups, let $G=G_1\times G_2$ be its direct product, and understand any subgroup of
$Aut(G_i)$ as a subgroup of $Aut(G)$ acting trivially on the other coordinate. It is easy to see that if the orbit
decidability for $A\leqslant Aut(G_1)$, and the conjugacy problem for $G_2$ are solvable, then the orbit decidability
for $A\leqslant Aut(G)$ is also solvable. And similarly, if $B\leqslant Aut(G_2)$ is orbit undecidable then so is
$B\leqslant Aut(G)$. Hence, answering question (a) in the affirmative reduces to answer in the affirmative questions
(b) and (c).

If $A$ is finitely presented, has trivial center, and CP($A$) is solvable, then we can take $G=A$, and the copy of $A$
in $Aut(G)$ given by conjugations is clearly orbit decidable. This answers (b) in the affirmative in this very
particular case.

Finally, let $B$ be a finitely generated, recursively presented group. By Higman's embedding theorem (see~\cite{LS}),
$B$ embeds in a finitely presented group $B'\neq 1$, which then embeds in $G_1 =B'\ast \mathbb{Z}$. Since $G_1$ has
trivial center, $Aut(G_1)$ contains a copy of $B$ given by inner automorphisms, say $B_1\leqslant Aut(G_1)$. Take now
another finitely presented group $G_2$ with unsolvable conjugacy problem, and consider their free product, $G=G_1 \ast
G_2$. Extend the morphisms in $B_1$ to automorphisms of $G$ acting trivially on $G_2$; this way, we obtain
$B_2\leqslant Aut(G)$, again isomorphic to $B$. Now, two given elements $u,v\in G_2$ lie in the same $(B_2\cdot
Inn(G))$-orbit if and only if they are conjugate to each other in $G_2$; thus, solvability of OD($B_2$) would imply
solvability of CP($G_2$). Hence, $B\simeq B_2 \leqslant Aut(G)$ is orbit undecidable. But, unfortunately, this does not
solve question (c) because, by construction, CP($G$) is unsolvable, like CP($G_2$).

Additionally, note that the recursive presentability for $B$ in the previous paragraph, is an extra condition also
satisfied in the main source of orbit undecidability presented above. Namely, all orbit undecidable subgroups coming
from Theorem~\ref{f2xf2} are of Mikhailova's type and so recursively presented (since they have solvable word problem).
At the time of writing we are not aware of any construction producing orbit undecidable subgroups which are not
recursively presented.  \qed

\section*{Acknowledgments}

We thank S. Hermiller, I. Kapovich and G. Levitt for interesting comments on the subject. The first named author is
partially supported by the INTAS grant N 03-51-3663 and by the grant Complex integration projects of SB RAS N 1.9. The
second and third named authors gratefully acknowledge partial support from the MEC (Spain) and the EFRD (EC) through
projects BFM2003-06613 and MTM2006-13544. The third one also thanks the Department of Mathematics of the University of
Nebraska-Lincoln for its hospitality during the second semester of the course 2003-2004, while a preliminary version of
this research was started. The three of them thanks the Centre de Recerca Matem\`atica (Barcelona, Catalonia) for its
hospitality during several periods of time while this research was conducted and written.

\end{document}